%% file: main.tex
\documentclass{article}
\usepackage[english]{babel}
\usepackage[belowskip=-5pt,aboveskip=2pt]{caption}
\usepackage{mathrsfs}
\usepackage{setspace}
 
\usepackage{arxiv}
\usepackage{parskip}
\usepackage{floatrow}
\usepackage[utf8]{inputenc} 
\usepackage[T1]{fontenc}    
\usepackage{url}            
\usepackage{booktabs}       
\usepackage{amsfonts}       
\usepackage{nicefrac}       
\usepackage{microtype}      
\usepackage{lipsum}
\usepackage{float}
\usepackage{subcaption}
\usepackage[toc,page]{appendix}

\usepackage[font=small,labelfont=bf,tableposition=top]{caption}
\usepackage{graphicx}
\usepackage{amsmath}
\usepackage[affil-it]{authblk}
\usepackage[belowskip=-5pt,aboveskip=2pt]{caption}
\usepackage{bm}
\usepackage{comment}
\usepackage{subcaption}
\usepackage{amsfonts}

\usepackage{abstract}

\usepackage{titlesec}  
\usepackage{nomencl}
\usepackage{titling} 
\usepackage [autostyle]{csquotes}
\usepackage{algorithm}
\usepackage{algpseudocode}
\usepackage{tikz}
\tikzset{
block/.style={
  draw, 
  rectangle, 
  minimum height=3cm, 
  minimum width=1cm,align=center,
  text width=8cm,
  execute at begin node=\setlength{\baselineskip}{20pt}
  }, 
line/.style={->,>=latex'}
}
\tikzstyle{myarrows}=[line width=10mm,draw=blue!60,-triangle 90,postaction={draw, line width=6mm, shorten >=5mm, -}]
\tikzset{line/.style={draw, thick, -latex'}}
\tikzset{every node/.style={font=\scriptsize}}
\usetikzlibrary{shapes.geometric}
\usetikzlibrary{arrows.meta}
\usetikzlibrary{tikzmark,calc}
\usetikzlibrary{shapes,arrows}
\usetikzlibrary{spy}
\usetikzlibrary{positioning}
\usetikzlibrary{backgrounds}
\usetikzlibrary{decorations}
\tikzset{
    custarr/.style={
        decorate, decoration={name=newarrow}%
    }
}

\usepackage[colorlinks=true,citecolor=blue, linkcolor=blue, urlcolor=black]{hyperref}
\AtBeginDocument{\hypersetup{pdfborder={0 0 1}}}
\usepackage{epstopdf}

\begin{document}

\title{ An iterative multi-fidelity approach for  model order reduction of multi-dimensional input parametric PDE systems}

\author{Manisha Chetry\textsuperscript{1*}, Domenico Borzacchiello\textsuperscript{1}, Lucas Lestandi \textsuperscript{1}, Luisa Rocha Da Silva\textsuperscript{1} }

\affil{\textsuperscript{1}  Nantes Université, Ecole Centrale Nantes, CNRS, GeM, UMR 6183, F-44000 Nantes, France} 
\affil{*\textbf{Corresponding author}:manisha.chetry@ec-nantes.fr}
\date{}

\renewcommand\Affilfont{\itshape\small}
\renewcommand{\maketitlehookd}{%

\begin{abstract}


We propose a parametric sampling strategy for the reduction of large-scale PDE systems with multidimensional input parametric spaces by leveraging models of different fidelity. The design of this methodology allows a user to adaptively sample points ad hoc from a discrete training set with no prior requirement of error estimators. It is achieved by exploiting low-fidelity models throughout the parametric space to sample points using an efficient sampling strategy, and at the sampled parametric points, high-fidelity models are evaluated to recover the reduced basis functions. The low-fidelity models are then adapted with the reduced order models ( ROMs) built by projection onto the subspace spanned by the recovered basis functions. The process continues until the low-fidelity model can  represent the high-fidelity model adequately for all the parameters in the parametric space. Since the proposed methodology leverages the use of low-fidelity models to assimilate the solution database, it significantly reduces the computational cost in the offline stage. The highlight of this article is to present the construction of the initial low-fidelity model, and a sampling strategy based on the discrete empirical interpolation method (DEIM).  We test this approach on a 2D  steady-state heat conduction problem for two different input parameters and make a qualitative comparison with the classical greedy reduced basis method (RBM), and further test on a  9-dimensional parametric non-coercive elliptic problem and analyze the computational performance based on different tuning of greedy selection of points.  
\end{abstract}

\keywords{ Multi-fidelity modeling, low-fidelity models, high-fidelity models, reduced basis method, greedy sampling, DEIM}
}
\maketitle

\section{\textbf{Introduction}}
\subsection{Motivation and background}

The convergence and efficiency of a reduced order model for approximation of the solutions of a large-scale PDE system depend heavily on the choice of the elements that constitute the "reduced basis" \cite{Priori}. Therefore, the input parameter selection for which the snapshots are generated must be appropriate. It must be sufficiently rich to adequately cover key areas of the parametric space. For sampling the parametric space, discretization techniques like uniform sampling and random sampling are frequently utilized. Both sampling methods, however, have their limitations. For instance, a multidimensional parametric system would necessitate a thorough grid search in uniform sampling, while random sampling would miss some crucial parts of the function in the parametric space. On the other hand, another sampling technique, known as Latin Hypercube Sampling (LHS)  \cite{helton2003latin,helton2005comparison} provides a compromise between uniformity and size of the sample, which makes it more efficient than uniform sampling and gives often better accuracy than random sampling. Some statistically-based sampling methods like  Monte Carlo methods are also among other popularly used techniques \cite{gentle2006random}. Even with efficient sampling techniques, the complexity can grow exponentially with the increase in the dimensionality of the parametric space, which is the case for many engineering applications. Therefore, the Reduced Basis Modeling (RBM) based on greedy sampling emerged as a promising tool for reducing the computational cost of Full Order Model (FOM) by generating high-fidelity snapshots at only a select few optimal parametric points. 

The basic idea in RBM is to adaptively choose sample points by finding the location at which the estimated error of the reduced model is maximum in the offline phase, thanks to rigorous error estimators. The greedy RBM was a subject of research for a very long time, first applied to find reduced models for the parameterized steady incompressible Navier– Stokes equations \cite{ito1998reduced} and then further developed for a variety of parameterized parabolic PDEs \cite{grepl2005posteriori,grepl2005reduced,grepl2007efficient} and also applied to several optimal control and inverse problems \cite{dede2010reduced,nguyen2010reduced}. During the past two decades, RBM had a significant contribution to the development of rigorous error bounds for Stokes flow problems, with a special focus on the inf sup stability conditions that can be referred to in the articles \cite{rozza2013reduced,martini2015reduced,quarteroni2007numerical}.

Even though, the RBM methods are not completely void of bottlenecks; it requires sharp, and  rigorous error estimators that are problem specific for ensuring the reliability of the method. Additionally, the error is estimated over a discrete training set, which must be a good surrogate of the continuous parametric space. Inefficient greedy sampling could yet occur from this, particularly for high dimensional parametric PDEs. To mitigate this issue, the authors in the article \cite{Cohen} have performed the greedy algorithm on random training sets of small sizes in every iterative cycle instead of estimating the error over the entire training set. The authors have successfully demonstrated a 16-parametric dimension system for a diffusion equation problem. Wilcox et al. in their research work \cite{bui2008model}, solved a sequence of optimization problems on a parametric space which is not a discrete set but rather continuous, to find optimal points adaptively in a greedy manner using both error estimator or residual error indicator. The authors demonstrated the proposed methodology on a thermal problem for the design of a fin based on two input parameters Biot number and conduction coefficients for  11, and 21 parametric dimensions respectively. The papers \cite{BUITHANH2007880,1582499,Hoang2015AnEG} also provide references to quite a few literary works on goal-oriented sampling methods. Numerous other studies have examined the use of subspace angles to assess the model's sensitivity to parameter changes \cite{Bazaz2015,varona2017automatic,amsallem2010interpolation} or sensitivity analysis to sample adaptively from the parametric space \cite{bond2007piecewise}.


\subsection{Overview of the idea: Iterative multi-fidelity model order reduction }


In this work, we suggest a sampling strategy that uses a multi-fidelity modeling approach as an alternative to the conventional greedy sampling technique that is driven by error estimators. Multi-fidelity modeling employs models of varied accuracy to estimate the same output quantity at a considerably faster convergence rate than using a single model with a higher approximation quality. For instance, in optimization problems, an optimizer is constructed to supply the input design variables at each iteration, and the model then assesses the value of the related objective function, its corresponding gradients, and the value of the constraint. Typically, high-fidelity models are retained in the loop to establish accuracy and convergence guarantees on the low-fidelity models, which are used to determine the best design parameters while taking advantage of computing speedups.  Low-fidelity model estimates are cheaper model approximations that can be a coarse-discretized model, a ROM, or even a  simplified physics model that can approximate the same output as the high-fidelity model but with lower accuracy.  Multi-fidelity methods have been successfully applied in optimization
problems based on co-kriging models \cite{cokriging}, uncertainty analysis \cite{uncertainty}, Monte Carlo simulations  \cite{WANG2021100035,montecarlo} to name a few. Such methods have also been extended to the machine learning and Physics-Informed Neural Network (PINN) domain, which can be referred to in  \cite{PENWARDEN2022110844,GUO2022114378}. The authors in \cite{PENWARDEN2022110844} have presented the connection of fidelity of different accuracy with neural networks by manipulating the width and depth of the network architecture.
 A comprehensive review of the past works and recent advances in the area of multi-fidelity modeling can be found in the survey \cite{peherstorfer2018survey}.

  Recently, in the work of Kast et al. \cite{KAST2020112947}, a multi-fidelity setup is exploited in the context of reduced-order modeling by solving a nonlinear structural dynamic problem based on three input parametric dimensions. First, an appropriate sampling set is exploited in the parametric space by leveraging a collection of low-fidelity models, followed by multi-fidelity Gaussian process regression (GPR) for approximation of the reduced coefficients in the online stage, therefore allowing an efficient decoupling of offline-online strategy. Another work based on a multi-fidelity adaptation strategy can be found in the work \cite{PEHERSTORFER201521} where the authors combined data-driven models with projection-based ROM and adapted the ROM for any change in the input parameter by low-rank updates to the reduced operators without rebuilding the FE operators from scratch. This work is focused on addressing the complexities of cases when the underlying properties of the PDE system are not static, but undergo dynamic changes due to the change in the latent variables.  
  
  In this work, we conjunct the multi-fidelity method with physics-based reduced order modeling for deriving low-cost ROMs efficiently without the need for problem-specific error estimators. Our approach is based on the "learn and adapt" framework. In the first step, a low-fidelity model learns to sample points from a large input parametric space, and in the second step, the low-fidelity model improves  by adapting to the current ROM approximation and the procedure continues until the low-fidelity model is a good representation of a FOM. By keeping  high-fidelity solutions in the loop, not only accuracy and convergence is achieved, the prior requirement of error estimators or upper bounds is no longer served.  The details are explained in section \ref{sec:IMFM}. The goal of this work is to explore the parametric space efficiently and generate appropriate snapshots for a high dimensional parametric system irrespective of the problem definition, the underlying discretization techniques used for solving the problem such as the finite element method (FEM), or finite volume method (FVM), or for cases when posteriori error estimators are unavailable or difficult to obtain.

In fact, we evaluate the max norm error between the high-fidelity solution and reduced order solution at the computed parametric points until it establishes the acceptable accuracy, hence our sampling process is adhoc based on heuristics. The greedy selection  of points can also be tuned as per the user's requirements for the efficient performance of the algorithm, which will be reflected in the numerical examples. We first introduce the algorithm on a 2D heat conduction problem with 2 input parameters and make a qualitative comparison with the existing greedy RBM. Further, we apply the proposed methodology to an advection-diffusion reaction problem in a parametric space of 9 dimensions.

This work focuses on two main aspects: the first is the construction of an initial low-fidelity model explained in section \ref{sec:initLF}, and the second is the sampling strategy using the DEIM technique explained in \ref{sec:DEIM}. The idea of using DEIM in iterative multi-fidelity modeling is not related to the approximation of the nonlinear term, but only to the notion of greedy selection of "optimal" parametric points.

The article is organized as follows: Section 2 presents a general framework of an elliptic partial differential problem, followed by classical reduced basis construction. Section 3 describes the construction of the reduced basis using the proposed iterative multi-fidelity approach. Section 4, demonstrates the methodology through different numerical examples followed by  results and discussion.

\section{\textbf{General problem setting}}
 This section presents a general problem setting for the purpose of establishing the notations for future reference. The proposed sampling methodology is unrelated to any particular PDE definition, but in order to compare it to the examples presented in the numerical analysis section \ref{numericaltest1} and \ref{numericaltest2}, we adhere to a general elliptic parameterized PDE. Let  $\Omega$ be some bounded domain and $\mathcal{V}$ be an associated functional space to characterize the solution as a function of space. 
Denoting with $\bm{\mu} $ as the input parameter, which belongs to the parametric space $\mathcal{D}$,
the problem is to find solution $u(\bm{\mu})$ in some finite-dimensional discrete  space $\mathcal{V}^\mathcal{N} \subset \mathcal{V} $ where $\mathcal{N}=\operatorname{dim}(\mathcal{V}^\mathcal{N})$ such that
\begin{equation}
a\left(u,v ; \bm{\mu}\right)=f\left(v; \bm{\mu}\right) \quad \forall v \in \mathcal{V}^{\mathcal{N}}
\label{PDE}
\end{equation}

In this work, finite elements are used to obtain the discrete solution to \ref{PDE}, however, the proposed method is also applicable to other numerical discretization solvers. After spatial discretization, the FE solution of the field variable `$u$', can be approximated as:

\begin{equation}
u \approx  u_h(\bm{x};\bm{\mu})= \sum_{i=1}^\mathcal{N} \mathrm{N}_i(\bm{x})(\mathrm{u}_h(\bm{\mu}))_i
\label{eq:FE}
\end{equation}

where, $\mathrm{N}_i$ are the shape functions of choice and $(\mathrm{u}_h(\bm{\mu}))_i$  represent the scalar values of the field $u_h$ at discretization points $\bm{x} \in \mathbb{R}^{\mathcal{N}}$. In practice, these values are stored as a collection of high-fidelity solutions that we call "snapshots"  of distinct parameter values $\mu\in \Xi_{train} \subset \mathcal{D}$, of cardinality $|\Xi_{train}|=N$. This set of snapshots $\{\bm{\mathrm{u}}_h^k\}_{k=1}^{N}$ is generated in the offline stage  by solving the PDE equation \eqref{PDE} using a high-fidelity solver for varying choice of input parameters. It can vary from material parameters to geometrical, and shape parameters, to boundary conditions of the unknown field variable, and so on.

\subsection{Reduced basis construction}

The basic idea in the reduced order modeling approach is that the discrete solution space associated with the underlying PDE lies in a low dimensional subspace and is in general represented by a low dimensional smooth manifold \cite{rozza2007reduced}. The manifold comprises all solutions of the parametric problem belonging to the parametric space i.e.
\begin{equation}\mathcal{M}_h=\operatorname{span}\left\{\bm{u}_h(\bm{\mu}) \in \mathcal{V}^\mathcal{N}: \bm{\mu} \in \mathcal{D}\right\} \subset \mathcal{V} \end{equation}
We aim to exploit a low dimensional structure of this manifold, $\mathcal{V}^R  \subset \mathcal{V}^\mathcal{N}$ where $R=\operatorname{dim}(\mathcal{V}^R)<<\mathcal{N}$ by appropriately choosing a set of orthogonal basis functions $\{\bm{\phi}_1,\bm{\phi}_2,\ldots,\bm{\phi}_R\}$  that spans the subspace of the manifold $\mathcal{M}_h$, and can well represent the manifold with small error. The associated reduced subspace is then given by,

$$\mathcal{V}^R=span\{\bm{\phi}_1,\bm{\phi}_2,\ldots,\bm{\phi}_R\} \subset \mathcal{V}^{\mathcal{N}}$$

The reduced solution $u^{R}\in \mathcal{V}^R $  can then be approximated by the linear combination of these basis functions given by:

\begin{equation}
u^{R}(\bm{x} ; \bm{\mu})=\sum_{i=1}^{R} \phi_i(\bm{x}) b^{i}( \bm{\mu})
\label{eq:RB}
\end{equation}

whose coefficients are calculated thanks to a projection  onto the reduced basis (RB) space. Therefore, the reduced problem can be sought as: for any $\bm{\mu}\in \mathcal{D}$, find $u^{R}(\bm{\mu}) \in \mathcal{V}^R$ such that:
\begin{equation}
 a\left(u^{R}, v_h; \bm{\mu}\right)=f(v_h ; \bm{\mu}) \quad \forall v_h \in \mathcal{V}^R
 \label{eq:ROM}
    \end{equation}

There are several strategies in the literature for constructing reduced basis functions, including the proper orthogonal decomposition (POD) \cite{sirovich1987turbulence,lumley1967structure} and the classical greedy reduced basis method. The objective of the greedy algorithm in the context of RBM is to adaptively enrich the reduced subspace with orthogonal basis functions \cite{Yvongreedy}. By doing so, the evaluation of high-fidelity snapshots for all the training parameters (as done in the classical POD) can be avoided in the offline step,  therefore reducing enormously the offline cost while improving the efficiency of MOR. It is based on the idea to select the parameter representing a local optimum in relation to an opportune error indicator iteratively,

\begin{equation}
\bm{\mu}_{n+1}=\operatorname{arg} \underset{\bm{\mu} \in \Xi_{train}}{ \operatorname{max}}  \Delta^R(\mu)
\label{greedyeq}
    \end{equation}

which means in the $(n+1)^{th}$ step, basically the sample point whose error metric $\Delta^R(\bm{\mu})$  indicates to be worst approximated among all the parameters $\bm{\mu} \in \Xi_{train} \subset \mathcal{D}$  by the solution of the current reduced model $\mathcal{V}_{n+1}^R$ is selected as the next sample point. At the sampled point, the high-fidelity snapshot is generated using the finite element approximation, followed by enrichment of the reduced basis subspace by Gram-Schmidt orthogonalization of the generated snapshots. This is repeated until the error estimator reaches a prescribed tolerance. To evaluate $\Delta^R(\bm{\mu})$ we need two essential ingredients: the dual norm of the residual and  a sharp lower bound for the coercivity constant which can be obtained by theta methods for simple PDEs (\cite{Hesthaven2015}) or successive constraint method (SCM) for general PDEs \cite{chen2009improved,chen2016certified,huynh2007successive}. 

It is essential for a good posterior error estimator to be sharp, and rigorous for ensuring the reliability of the RBM and also has to be  computationally inexpensive for efficient greedy sampling. However, for complex PDEs  construction of sharp and rigorous error bounds may not be easily achievable which can lead to unreliable reduced basis approximation, and also for PDE systems with high dimensional parametric space, finding an error estimate over the entire parametric set may become computationally expensive. As a result, a novel technique is suggested in this work, where we use POD for the basis construction  but a different strategy utilizing multi-fidelity modeling is employed to effectively choose the snapshots, as explained in the next section. 


\section{Iterative multi-fidelity modeling (IMF) for building PODG-ROM}\label{sec:IMFM}

Traditionally, if $\{\bm{\mathrm{u}}_{\text{\tiny HF}}^k(\bm{\mu}) \}_{k=1}^N \in \mathbb{R}^{\mathcal{N}}$ represents the high-fidelity snapshots of the PDE problem \ref{PDE} at distinct parameter values $\bm{\mu}\in \Xi_{train} \subset \mathcal{D}$, of cardinality $|\Xi_{train}|=N$,  that can be suitably well approximated in a low-dimensional manifold, then the solution can be represented in a separated form as,

\begin{equation}
u_{\text{\tiny HF}}(\bm{x};\bm{\mu}) \approx \sum_{i=1}^{r}{\phi}_{\text{\tiny HF}}^i(\bm{x})\psi_{\text{\tiny LF}}^i (\bm{\mu}) 
\label{eq:HF}
\end{equation}

Here, $\{\bm{\phi}_{\text{\tiny HF}}^i\}_{i=1}^r \in \mathbb{R}^\mathcal{N}$ represents the high-fidelity basis functions  
which spans the low dimensional subspace, i.e. $\mathcal{V}^R= \operatorname{span}\{\bm{\phi}_{\text{\tiny HF}}^1, \dots, \bm{\phi}_{\text{\tiny HF}}^r\} \subset \mathcal{V}^{\mathcal{N}}$, $R= dim(\mathcal{V}^R)$   and $\psi_{\text{\tiny LF}}(\bm{\mu}):\mathcal{D}\longrightarrow \mathbb{R}$ are parametric functions that span the parametric space. It is true that high-fidelity models can capture the intricacies of complex PDE systems, but they are also equally expensive to train and the offline cost to recover the basis functions is quite high. If the parametric functions $\psi_{\text{\tiny LF}}^i(\bm{\mu})$ were previously known, we could easily extract a set of points $\bm{\mu}^{\bm{P}} \subset \Xi_{train}$ where  $\bm{P} \in (1,N)$ using any efficient sampling technique. The optimal basis functions might then be recovered by generating high-fidelity snapshots at the computed set of points. 

However, $\psi_{\text{\tiny LF}}^i (\bm{\mu})$  are not known a priori, but we can reasonably assure that if a low-fidelity model is instead used for  approximation, then by similar expression as given in eq.\ref{eq:HF} we have,

 
\begin{equation}
u_{\text{\tiny LF}}(\bm{x};\bm{\mu}) \approx \sum_{i=1}^{r}{\phi}_{\text{\tiny LF}}^i(\bm{x})\psi_{\text{\tiny LF}}^i (\bm{\mu}) )
\label{eq:LF}
\end{equation}

 where $\psi_{\text{\tiny LF}}^i (\bm{\mu})$ have similar features as the $\psi_{\text{\tiny HF}}^i(\bm{\mu})$. Therefore, we can use those to recover the high-fidelity snapshots at computed points  $u_{\text{\tiny HF}}(\bm{\mu}^{\bm{P}})$. 
The process can be made iterative, as the newly computed $\bm{\mu}^{\bm{P}}$ can effectively result in the reconstruction of the high-fidelity basis functions. This subsequently leads to the enrichment of the reduced basis subspace, which causes an improvement in the low-fidelity model approximation.
 
\begin{figure}[h!]
\centering
\resizebox{0.8\columnwidth}{!}{
\tikzstyle{every node}=[font=\fontsize{42}{42}\selectfont]
     \tikzstyle{block} = [rectangle split, draw, rectangle split parts=2,text width=60em, rounded corners, minimum width=40em, minimum height=20em, align=center,font=\fontsize{42}{42}\selectfont]

    \tikzstyle{blueblock} = [rectangle, draw, fill=blue!20, minimum width=40em, text width=60em, text centered, rounded corners, minimum height=20em, align=center,font = \fontsize{42}{42}\selectfont] 
    
    \tikzstyle{whtblock} = [rectangle,rounded corners, draw, fill=white!20, text width=24em, minimum height=10em, align=left, font = \fontsize{40}{40}\selectfont]  
\tikzstyle{decision} = [diamond, draw]

 \begin{tikzpicture}[fill=gray, scale = 5, transform shape]
\input{sequence}

\end{tikzpicture}
}

\vspace*{0.5cm}
\caption{Flow of the iterative multi-fidelity modeling approach.}
\label{fig:IMF}
\end{figure}
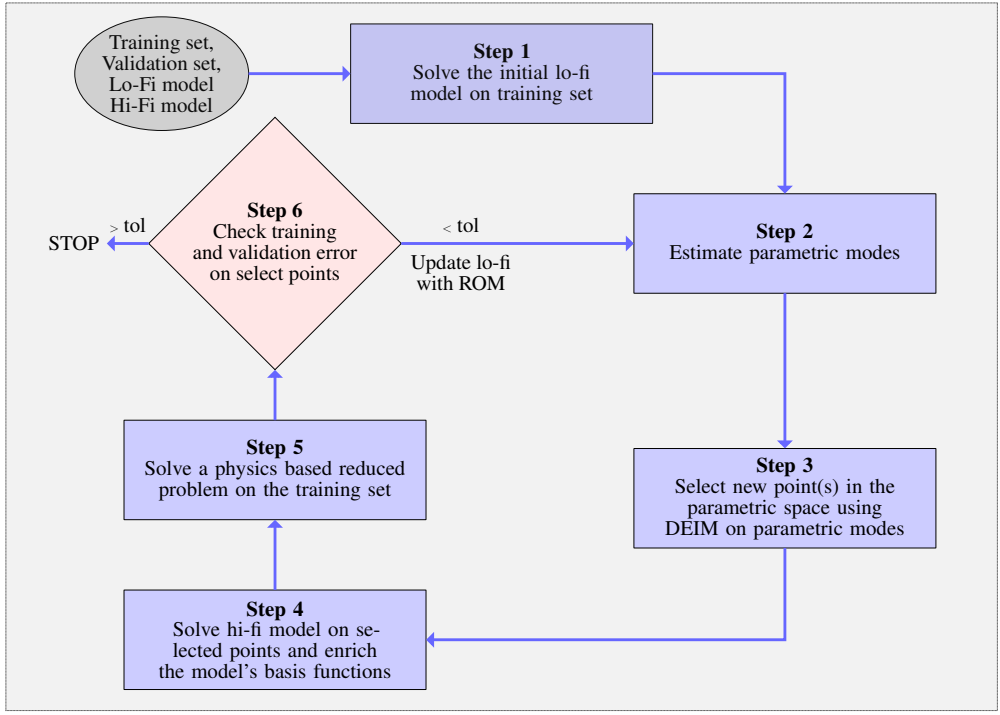

Therefore, step 1 of the proposed method is to obtain a poor or inexpensive approximation to the FOM using a low-fidelity (lo-fi) model, $f_{\text{\tiny LF}}^{init}:\Omega \times \mathcal{D} \rightarrow \mathbb{R} $  that maps all the parameters belonging to a given training set $\bm{\mu} \in \Xi_{train} \subset \mathcal{D}$ to produce the same output with lower accuracy. The flow is shown in figure \ref{fig:IMF} and details on the construction of the initial low-fidelity model are explained in section \ref{sec:initLF}.

In step 2, we compute parametric functions $\{\bm{\psi}_{\text{\tiny LF}}^i\}_{i=1}^r \in  \mathbb{R}^N$ from the low-fidelity model  approximation by primarily using POD, and then in step 3, we recover "optimal sampling points" using interpolation strategy over these parametric functions to sample optimal points in a greedy procedure, $\bm{\mu}^{\bm{P}} \subset \Xi_{train}$ with $\bm{P} \in (1,N)$. The details are discussed in section \ref{sec:DEIM}. 

Next, in step 4, using a high-fidelity (hi-fi) model we generate snapshots on the select sample points $\bm{\mu}^{\bm{P}}$ to recover the high-fidelity basis functions $\bm{\phi}_{\text{\tiny HF}}^i \in \mathcal{V}^R$  and therefore, enrich the reduced basis subspace,  $\mathcal{V}^{R} =\operatorname{span}(\{\bm{\phi}_{\text{\tiny HF}}^i\}_{i=1}^{r''}) \in \mathbb{R}^\mathcal{N}$ where $r''\leq \#(\bm{\mu}^{\bm{P}})$, the construction is properly explained in the section \ref{sec:basisfunc}. A high-fidelity model is mathematically defined as $f_{\text{\tiny HF}}:\Omega \times \mathcal{D} \rightarrow \mathbb{R} $ that maps all the selected points $\bm{\mu}^{\bm{P}} \subset  \Xi_{train} $ to estimate the output with the accuracy that is needed for the task. 

In step 5, we then solve a POD-G ROM for all $\bm{\mu} \in \Xi_{train}$. The functions $\psi_{\text{\tiny LF}}^i$ obtained during the first approximation are likely to not generate exactly the same space as the $\psi_{\text{\tiny HF}}^i$ i.e $\operatorname{span}(\psi_{\text{\tiny LF}}^i) \neq \operatorname{span}(\psi_{\text{\tiny HF}}^i) $. Therefore, recovery of the high-fidelity basis functions $\bm{\phi}_{\text{\tiny HF}}^i$ may not be accurate, and reliable to represent the large-scale PDE system, hence the procedure has to undergo certain iterations. 

 Finally, in step 6, we evaluate the error between the high-fidelity model and reduced basis model approximation  at the computed discrete points using the error metrics that are discussed in \ref{errordef}. If it is below a certain prescribed tolerance level,  we terminate the algorithm, else we adapt the low-fidelity model with the current ROM approximation (refer section \ref{sec:LFupdate}) and repeat the procedure until the ROM constructed represents the FOM adequately. To  measure the overall performance of the algorithm,  we also check for validation error by computing the error on another set of parameters belonging to a given validation set $\Xi_{val} \subset \mathcal{D}$.


The sample points that are obtained  provide a locally optimal choice at each stage of the iterative cycle, however as iteration  continues and new points are added in each iteration the algorithm converges towards the global solution with certain accuracy in very reduced complexity. 

The proposed method in this work is tested on an elliptic PDE problem. It can be extended to hyperbolic or parabolic PDEs as well. However, the reduced basis subspace needs to be built appropriately to take into account the time integration.

\subsection{Construction of the initial low-fidelity model}\label{sec:initLF}
We propose two types of sketch models for the construction of an initial low-fidelity model, depending on the availability of either of the two conditions:

\begin{enumerate}
  \item No database of high-fidelity solutions is present a priori, 
   \item A database of high-fidelity solutions is available

   \end{enumerate}

\paragraph{Coarse sketch model:} When no database of solutions is present a priori, the initial low-fidelity model is built using a  derefined version of a high-fidelity model, which is nothing but a coarse finite element model. We refer to it as a "coarse sketch model" for future discussions. Snapshots generated using this coarse sketch appear to be a wide matrix, denoted by $\bm{S}_{\text{ \tiny LF}}=[\bm{\mathrm{u}}_{\text{\tiny LF}}^1,\bm{\mathrm{u}}_{\text{\tiny LF}}^2,\ldots, \bm{\mathrm{u}}_{\text{\tiny LF}}^N] \in \mathbb{R}^{m\times N}$, where $\bm{\mathrm{u}}_{\text{\tiny LF}}=\bm{\mathrm{u}}_h(\bm{\mu})$ at discretization points $\bm{x}_{\text{\tiny LF}} \in \mathbb{R}^m$, such that $ m << \mathcal{N}$. 

\paragraph{POD-G ROM:} Now, let us assume the case when we already have some solutions of the  large-scale PDE system available to us, be it experimental or numerical data.  This is materialized in our work by solving a high-fidelity model for any random training parameters, $\mathbf{X}=[\bm{\mathrm{u}}_{\text{\tiny HF}}^1(\bm{\mu}^k),\bm{\mathrm{u}}_{\text{\tiny HF}}^2(\bm{\mu}^k),\ldots, \bm{\mathrm{u}}_{\text{\tiny HF}}^K(\bm{\mu}^k)] \in  \mathbb{R}^{\mathcal N \times K}$ with $ \mathcal{N} >>K$, and $k\in (1,N)$ and term it here as a \textbf{random sketch model}. Then the initial low-fidelity model is a ROM approximation, which is constructed by Galerkin projection of the PDE system onto the reduced basis functions computed from this sketch model. 

\paragraph{Remark: }  Although the selection parameters for a random sketch model are entirely user-dependent, it is advised to start by building the reduced bases with just a few snapshots. For the primary purpose of lowering the offline cost of MOR, the low-fidelity model approximation from such a random sketch must remain a less expensive approximation to the FOM at the initial stage of the method. Then, as the iteration advances, the quality of the low-fidelity model improves and converges to the FOM accurately.

 \subsection{Parametric point selection}\label{sec:DEIM}

 As previously mentioned, sampling points are extracted from the parametric modes of the low-fidelity model approximation. This is based on the heuristic assumptions that the low-fidelity model, although a poor approximation to the high-fidelity model, may nonetheless accurately reflect the essential features of the high-fidelity model's parametric dependence. A brief description of how the parametric points are sampled from a given training set, $\bm{\mu}^{\bm{P}} \subset \Xi_{train} $ using DEIM is explained in the algorithm \ref{algo:deim}. DEIM finds the sample points in a greedy way from an input basis which is given here by the parametric functions computed by performing SVD on the low-fidelity model approximation,
    
 \begin{equation}
 \bm{S_{\text{\tiny LF}}}=\bm{\Phi}_{\text{\tiny LF}} \bm{\Sigma}_{\text{\tiny LF}} \bm{\Psi}_{\text{\tiny LF}}^{\text{\tiny T}}
 \end{equation}
 
  where, $\bm{\Phi}_{\text{\tiny LF }}=[\bm{\phi}_{\text{\tiny LF }}^1,\bm{\phi}_{\text{\tiny LF }}^2,\ldots,\bm{\phi}_{\text{\tiny LF }}^r] \in \mathbb{R}^{m \times r}$ with $m<<\mathcal{N}$, and $r \leq N$ are poorly approximated POD modes while $ \bm{\Psi}_{\text{\tiny LF}}=[\bm{\psi}_{\text{\tiny LF}}^1,\bm{\psi}_{\text{\tiny LF}}^2,\ldots,\bm{\psi}_{\text{\tiny LF}}^r] \in \mathbb{R}^{N \times r}$ denotes the parametric modes that span the parametric subspace. The rectangular diagonal matrix $\bm{\Sigma}_{\text{\tiny LF}} \in \mathbb{R}^{r \times r}$ contains the corresponding non-negative singular values, $\sigma_1 \geq \ldots \geq \sigma_r \geq 0$ accounting to the information content of the low-fidelity model solution data. The process of sampling starts by selecting the index  with the largest magnitude, corresponding to the first entry of the input basis $\{\bm{\psi}_{\text{\tiny LF}}^i\}_{i=1}^r$. The remaining points are selected by finding the location at which the residual of the current approximation is maximum (refer to the algorithm \ref{algo:deim}). The points that are computed are unique due to the linear independence of the input basis, which guarantees that the indices are hierarchical and non-repetitive in nature. 

 \begin{algorithm}
\caption{DEIM sampling adopted from \cite{Chaturantabut2010}}
\begin{algorithmic}[1]
\State INPUT: Parametric functions $\bm{\Psi}_{\text{\tiny LF}} \in \mathbb{R}^{N \times r}$  
\State OUTPUT: Sampled parametric points 
 $\bm{\mu}^{\bm{P}}=\left[\bm{\mu}^{(1)}, \ldots, \bm{\mu}^{(r)}\right]^{\text{\tiny T}} \subset \Xi_{train}$

\State  $ \quad \bm{\mu}^{(1)}=\max \left\{\left|\bm{\psi}_{\text{\tiny LF}(\cdot,1)}\right|\right\}$
\State $ \quad \bm{\mu}^{\bm{P}}=[\bm{\mu}^{(1)}]$

\For{ $l=2:r$ } 

\State $\quad \quad \text{Solve }c= [\bm{\Psi}_{\text{\tiny LF}(1:l-1,1:l-1 )}]^{-1} \bm{\psi}_{\text{\tiny LF}(1:l-1,l)}$ 

\State $\quad \quad \bm{\mathrm{r}}=\bm{\psi}_{\text{\tiny LF}(\cdot,l)}-\bm{\Psi}_{\text{\tiny LF}} c$

\State $\quad \quad \bm{\mu}^{(l)}=\max \{|\bm{\mathrm{r}}|\}$

\State $\qquad \bm{\mu}^{\bm{P}}=[\bm{\mu}^{\bm{P}},\mu^{(l)}]^{\text{\tiny T}}$

\EndFor
\end{algorithmic}
\label{algo:deim}
\end{algorithm}

This sampling procedure is resumed after every iteration of the proposed algorithm and is not restarted from the beginning. As a result, instead of oversampling the same points from the training set, we are able to sample distinct points from it. The parametric functions obtained in $i^{th}$ iteration of the multi-fidelity algorithm are orthogonalized with respect to the parametric functions obtained in $(i-1)^{th}$ iteration through Gram-Schmidt orthonormalization in order to prevent repetition and picking up points closer to previously calculated points. This step is crucial for finding the best and most distinctive points throughout each iteration cycle, enabling us to explore the parametric space more thoroughly.
  \paragraph{Remark: }A general remark to take into account that while picking parametric points, it is always a better choice to select from the first `$r$' truncated parametric functions due to its content of the highest energy or information of the system arranged in descending order.
  
  It is noteworthy that the proposed methodology doesn't necessarily perform like a classical greedy sampling procedure due to the nature of the selection of points from the parametric functions. Hence, depending on the available computing resources, the level of "greediness" can be fine-tuned. In other words, since DEIM will generate the same number of points as the rank of parametric functions, the user can decide to select all the sample points at once, or also has the option to select one parametric point per iteration. 
 This is one of the advantages of this approach where the selection of points per iteration is completely user-dependent, which can be assimilated for parallel computations. A general remark has to be made, incorporating such a step can also lead to the selection of excess sample points than required which deviates from the main objective of selecting a few optimal points and hence needs to be taken care of.

 \subsection{Recovery of the reduced basis functions}\label{sec:basisfunc}
The high-fidelity basis function in the first iteration of the proposed method is recovered by performing SVD on the select snapshots, $\bm{S}_{\text{\tiny HF}}=\{\bm{\mathrm{u}}_{\text{\tiny HF}}(\bm{\mu}^{\bm{P}})\}$  for all $\bm{\mu}^{\bm{P}} \subset \Xi_{train}$.

 \begin{equation}
\operatorname{svd}(\bm{S_{\text{\tiny HF}}})= \bm{\Phi}_{\text{\tiny HF}} \bm{\Sigma}_{\text{\tiny HF}} \bm{\Psi}_{\text{\tiny HF}}^{\text{\tiny T}}
 \end{equation}

where, $\bm{\Phi}_{\text{\tiny HF}}=[\bm{\phi}_{\text{\tiny HF}}^1,\bm{\phi}_{\text{\tiny HF}}^2,\ldots,\bm{\phi}_{\text{\tiny HF}}^{r''}] \in \mathbb{R}^{\mathcal{N}\times r''}$ with  $ r''\leq \#(\bm{\mu}^{\bm{P}})$ contains the high-fidelity reduced bases that span the low-dimensional subspace $\mathcal{V}^{R}$ and $ \bm{\Psi}_{\text{\tiny HF}}=[\bm{\psi}_{\text{\tiny HF}}^1,\bm{\psi}_{\text{\tiny HF}}^2,\ldots,\bm{\psi}_{\text{\tiny HF}}^{r''}] \in \mathbb{R}^{N \times r''}$ denotes the parametric modes. Similarly, the rectangular diagonal matrix $\bm{\Sigma}_{\text{\tiny HF}} \in \mathbb{R}^{r'' \times r''}$ contains the corresponding non-negative singular values $\sigma_1 \geq \ldots \geq \sigma_{r''} \geq 0$ accounting to the information content of the high-fidelity model solution data. At $(i+1)^{th}$ iteration of the algorithm, the reduced subspace $\mathcal{V}^{R}$ is updated through the Gram-Schmidt procedure (refer algorithm \ref{algo:GS}).

\begin{algorithm}
\caption{ Gram-Schmidt orthonormalization at the $(i+1)^{th}$ iteration of the proposed method}
\label{algo:GS}

\begin{algorithmic}[1]

\For{ $l=1:\text{dim}(\bm{S}_{\text{\tiny HF}})$}
\State   $\bm{\phi}_{\text{\tiny HF}}^l= \bm{S}_{\text{\tiny HF}}^l-\mathcal{V}_i^{R} \langle \mathcal{V}_i^{R},\bm{S}_{\text{\tiny HF}}^l \rangle $
\If{$\frac{||\bm{\phi}_{\text{\tiny HF}}^l||}{||\bm{S}_{\text{\tiny HF}}^l||} > \epsilon_g$}
\State  \qquad $\mathcal{V}_{i+1}^{R}= \mathcal{V}_i^{R} \bigoplus \frac{\bm{\phi}_{\text{\tiny HF}}^l}{||\bm{\phi}_{\text{\tiny HF}}^l||}$
\EndIf 
\State   Update, $\mathcal{V}_i^{R}=\mathcal{V}_{i+1}^{R}$
\EndFor
\end{algorithmic}
\end{algorithm}

\subsection{Updating low-fidelity model }\label{sec:LFupdate}

In this part, we demonstrate how the current POD-G ROM approximation can be used to update the low-fidelity model for each iteration of the algorithm until convergence, i.e. we approximate the solution $u_{\text{\tiny LF}}:\Omega \times \mathcal{D} \rightarrow \mathbb{R}$ with a function $u^R\in \mathcal{V}^R$  defined by,

\begin{equation}
u^R(\bm{x};\bm{\mu} )=\sum_{i=1}^{r''}\bm{\phi}_{\text{\tiny HF}}^i (b(\bm{\mu}))_i= \bm{\Phi}_{\text{\tiny HF}}\bm{b}(\bm{\mu})\\ 
 \label{eq:HFRB}
\end{equation}

where, the POD expansion coefficients $\bm{b}(\bm{\mu})=(b_1,b_2,\ldots,b_{r''})^{\text{\tiny T}}$ can be calculated by Galerkin projection of the PDE system onto the basis functions $\bm{\phi}_{\text{\tiny HF}}^i$. The initial low-fidelity snapshots data is now updated with the current reduced solution, such that $\bm{S}_{\text{\tiny LF}}=[u^R_1,u^R_2,\ldots, u^R_N] \in \mathbb{R}^{\mathcal{N}\times r''}$. 


\textbf{Remark}: To improve the efficiency of the method, one can also approximate the low-fidelity data with the coefficients of the POD expansion instead of the reduced solution itself and replace $\bm{S}_{\text{\tiny LF}}$ with $\bm{B}$ where $\bm{B}=[\bm{b}^1,\bm{b}^2,\ldots, \bm{b}^N] \in \mathbb{R}^{{r''}\times N}$. Here $\bm{b}^i=(b_1,b_2,\ldots,b_{r''})^{\text{\tiny T}}$ represents the POD coefficients from eq. \ref{eq:HFRB}.  This process can reduce the cost of exploration of the parametric space using low-fidelity approximation from $\mathcal{O}(\mathcal{N})$ to $\mathcal{O}(r'')$.

By performing SVD on $\bm{B}$ we have,

\begin{equation}
\bm{B}= \bm{\varphi}  \bm{\varsigma}\hat{\bm{\psi}}^{\text{\tiny T}}
 \label{eq:svdLF}
 \end{equation}

Rewriting eq. \ref{eq:HFRB} we have,
\begin{equation}
\tilde{u}^{R}(\bm{x};\bm{\mu} )\approx u^{R}(\bm{x};\bm{\mu} )=\underset{\bm{\tilde{\Phi}}_{\text{\tiny HF}}}{\underbrace{\bm{\Phi}_{\text{\tiny HF}} \bm{\varphi}}}\bm{\varsigma}\hat{\bm{\psi}}^{\text{\tiny T}}
 \label{eq:LFRB}
\end{equation}

The original high-fidelity basis functions can now be replaced by the  approximate left singular vectors by the expression shown in eq. \ref{eq:LFRB}. This scaling factor is committed to improving the accuracy of the basis functions that could be lost if reduced coefficients are used in place of reduced-order solutions for the parametric exploration. 

The next section is dedicated to  a discussion on error metrics which are used to test the reliability of the approach.
 
\subsection{Error metric}\label{errordef}
If $u_{\text{\tiny HF}}(\bm{x};\bm{\mu}^{\bm{P}})$  and  $u^{R}(\bm{x};\bm{\mu}^{\bm{P}})$ represents the FOM and ROM solution respectively at the computed parametric points $\bm{\mu}^{\bm{P}}$ from the training set $\Xi_{train} \subset \mathcal{D}$, then the max norm of the relative error is estimated at the sample points such that 
\begin{equation} \epsilon_{train}=\underset{\bm{\mu}^{\bm{P}}}{max} \sqrt{  \frac{\sum_{ i = 1}^{\mathcal{N}}||{u_\text{\tiny HF}(x_i;\bm{\mu}^{\bm{P}})-u^{R}(x_i;\bm{\mu}^{\bm{P}})||}_2^2}{\sum_{ i = 1}^{\mathcal{N}} ||u_{\text{\tiny HF}}(x_i;\bm{\mu}^{\bm{P}})||_2^2}}  \end{equation}

In order to better understand the quality of the reduced model, the algorithm is validated on another set of parameters, $ \Xi_{val} \subset \mathcal{D} $ for the problem defined. By validating on a different set of points, if the error estimation between the FOM and ROM solutions is reduced as the cycle increases, it can be inferred that ROM well approximates the large-scale system for any $\bm{\mu} \in \mathcal{D}$. 
If the max norm error for both the training and validation set is below a certain tolerance limit, the sketch model constructed can be considered to be reliable. 

\begin{equation} \epsilon_{val}= \underset{\bm{\mu} }{max} \sqrt{  \frac{\sum_{ i = 1}^{\mathcal{N}}||{u_{\text{\tiny HF}}(x_i;\bm{\mu} )-u^{R}(x_i;\bm{\mu} )||}_2^2}{\sum_{ i = 1}^{\mathcal{N}} ||u_{\text{\tiny HF}}(x_i;\bm{\mu} )||_2^2}} \qquad \forall \bm{\mu} \in \Xi_{val}\end{equation}

The error between all the snapshots and the ROM solutions obtained by the iterative multi-fidelity approach is evaluated  and can be viewed as a benchmark for the ROM error, provided sufficient snapshots are generated:
\begin{equation} \epsilon_{\text{\tiny ROM}}= \sqrt{\sum_{ i = 1}^{N}  \frac{||{u_{\text{\tiny HF}}(\bm{x};\bm{\mu}^{(i)})-u^{R}(\bm{x};\bm{\mu}^{(i)})||}_2^2}{ ||u_{\text{\tiny HF}}(\bm{x};\bm{\mu}^{(i)})||_2^2}} \qquad  \bm{\mu} \in \Xi_{train}\end{equation}

We also analyze the POD basis projection error, which is given by the error between the snapshots and their projection onto the recovered basis functions:

\begin{equation} \epsilon_{\text{\tiny POD}}= \sqrt{\sum_{ i = 1}^{N}  \frac{||{u_{\text{\tiny HF}}(\bm{x};\bm{\mu}^{(i)})-\Pi u_{\text{\tiny HF}}(\bm{x};\bm{\mu}^{(i)})||}_2^2}{ ||u_{\text{\tiny HF}}(\bm{x};\bm{\mu}^{(i)})||_2^2}} 
\qquad  \bm{\mu} \in  \Xi_{train}\end{equation}

\section{ 2D heat conduction problem}
\label{numericaltest1}
\begin{figure}[hbt!]
    
    \centering
       \includegraphics[width=0.5\textwidth]{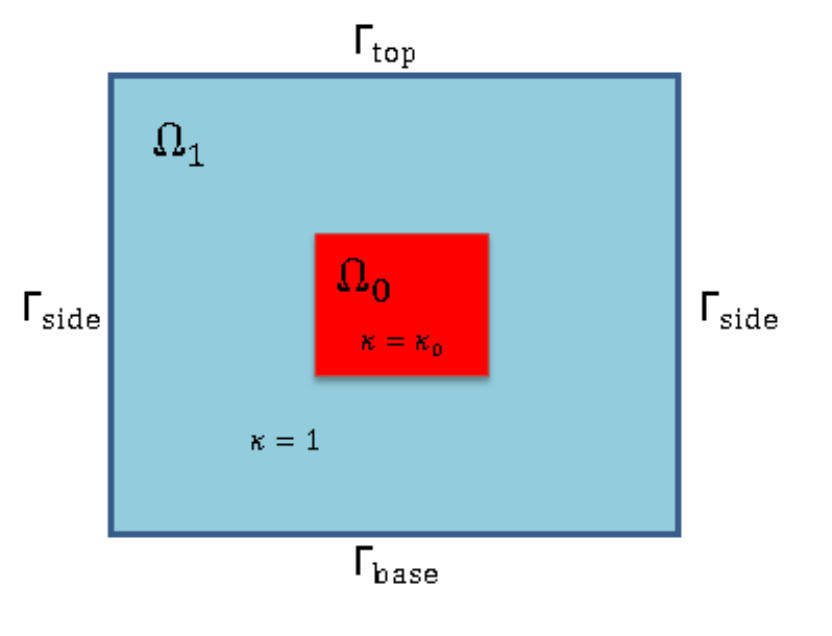}
       \caption{Geometry of heat conduction problem}
      \label{conduction}
  \end{figure}

In order to assess the proposed methodology, we first begin by analyzing a simple steady-state heat conduction  problem in a 2D domain,  $\Omega=(0,1) \times(0,1)$ as shown in figure \ref{conduction}. This problem is reproduced from (\cite{Hesthaven2015}) where it is solved using classical greedy RBM. The boundary of the domain is split into three parts, the base, the top, and the sides, and $\Omega_{0}$ is a square block placed in the center of the domain. Let $\kappa$ be the thermal conductivity with $\kappa|_{\Omega_{0}}=\mu_{[1]}$ and $\kappa|_{\Omega_{1}}=1$ where  $\Omega_{1}=\Omega \backslash \Omega_{0}.$ 

Two input parameters are considered for this problem $\bm{\mu}=[\mu_{[1]},\mu_{[2]}]$, where $\mu_{[1]}$ is the conductivity in $\Omega_0$ the region, and the second parameter $\mu_{[2]}$ is the constant heat flux over the bottom boundary. The strong formulation for this parameterized problem is governed by Poisson's equation. For some parameter value $\bm{\mu} \in \mathcal{D}$, find $u(\bm{\mu})$ such that

\begin{equation}
\begin{aligned}
\nabla \cdot \kappa_{\mu} \nabla u(\bm{\mu}) &=0 & & \text { in } \Omega, \\
u(\bm{\mu}) &=0 & & \text { on } \Gamma_{\text {top }} \\
\kappa_{\mu} \nabla u(\bm{\mu}) \cdot n &=0 & & \text { on } \Gamma_{\text {side}} \\
\kappa_{\mu} \nabla u(\bm{\mu}) \cdot n &=\mu_{[2]} & & \text { on } \Gamma_{\text {base}}
\end{aligned}
\end{equation}

\begin{figure}[hbt!] 
 \centering
    \begin{subfigure}{0.4\textwidth}
     \includegraphics[trim=6cm 11cm 6cm 11cm, clip=true,width=\textwidth]{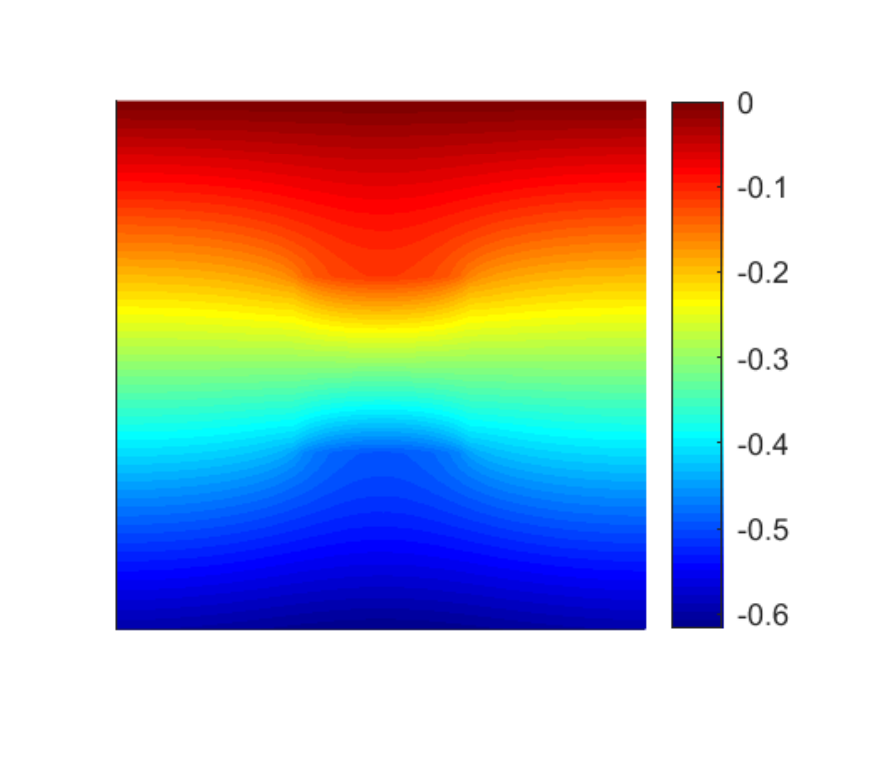}
     \caption{  $\mu_{[1]}=0.1, \mu_{[2]}=-1$}
  \end{subfigure}
 \begin{subfigure}{0.4\textwidth}
   \includegraphics[trim=6cm 11cm 6cm 11cm, clip=true,width=\textwidth]{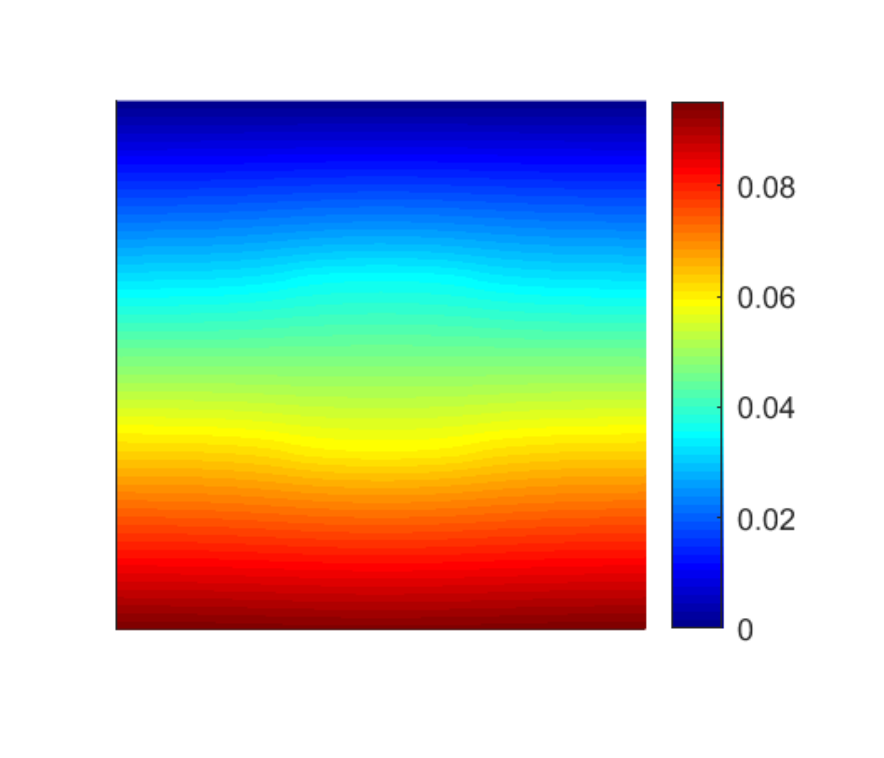}
   \caption{$\mu_{[1]}=0.4, \mu_{[2]}=0.5 $}
        
  \end{subfigure}
  
     \caption{Two different representative solutions for the parameterized conductivity problem.}
       \label{conductionprofile}
\end{figure}

Here,  $u(\bm{\mu})$  is the scalar  temperature field  variable, and $\kappa_{\mu}$ is given such that  $\kappa_{\mu}=\varphi_{1}+\mu_{[1]} \varphi_{0}$, where $\varphi$ is the characteristic function with subscript donating the corresponding domain.  
Defining $\mathcal{V}^\mathcal{N}=\left\{v \in H_0^1(\Omega)|v|_{\Gamma_{\text {top }}}=0\right\}$, the weak parametrized formulation then reads: for some parameter $\mu \in \mathcal{D}$, find $u(\mu) \in \mathcal{V}^\mathcal{N}$ such that, 

\begin{equation}
\begin{gathered}
a(u(\mu), v ; \bm{\mu})=f(v ; \bm{\mu}) \quad \forall v \in \mathcal{V}^\mathcal{N}, \\
a(w, v ; \bm{\mu})=\int_{\Omega} \kappa_{\bm{\mu}} \nabla w \cdot \nabla v \text { and } f(v ; \bm{\mu})=\mu_{[2]} \int_{\Gamma_{\text {base }}} v,
\end{gathered}
\label{eq:weakformconduction}
\end{equation}

for all $v, w \in \mathcal{V}^\mathcal{N}$. The selected range for the parametric study: $\bm{\mu}=[\mu_{[1]}, \mu_{[2]}]\in \mathcal{D}=[0.1,10] \times[-1,1]$. A total of, 2050 sample points are generated in which the training set $\Xi_{train}$ compromises of 2000 points and the validation set $\Xi_{val}$ consists of 50 points. For $\mu_{[1]}$, the points are generated using uniform discretization whereas for the second input parameter $\mu_{[2]}$, the points are generated using log space. The graphical representation of the temperature field  for two different sets of parameters is shown in figure \ref{conductionprofile}.

\subsection{Results and discussion}
\begin{figure}[hbt!] 
    \begin{subfigure}{0.4\textwidth}
     \includegraphics[trim=6cm 11cm 6cm 11cm, clip=true,width=\textwidth]{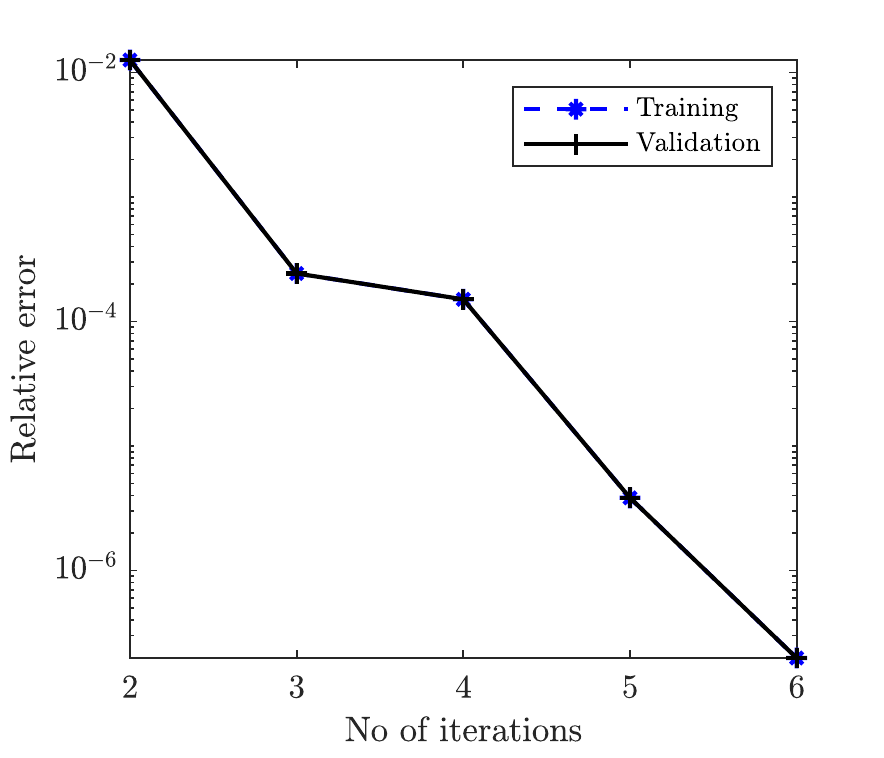}
     \subcaption{Initial random sketch}
          \label{err1a}

  \end{subfigure}
 \begin{subfigure}{0.4\textwidth}
   \includegraphics[trim=6cm 11cm 6cm 11cm, clip=true,width=\textwidth]{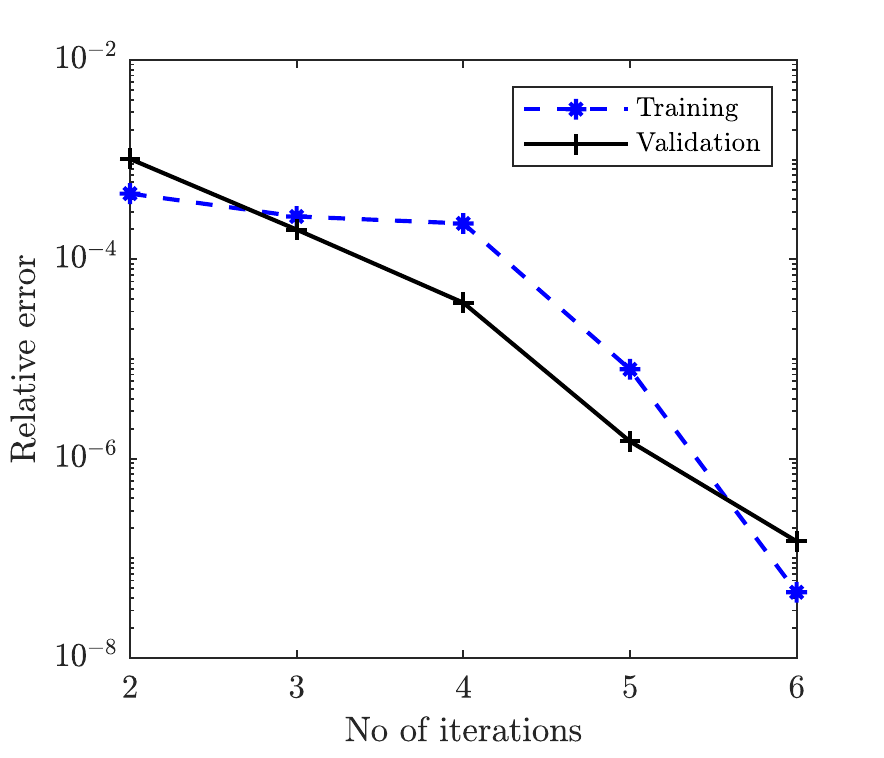}
   \subcaption{Initial coarse sketch}
   \label{err1b}

  \end{subfigure}

     \caption{ Error $\epsilon_{train}$ and $\epsilon_{val}$ between FOM and ROM solution using  random sketch model (\textbf{a}), coarse sketch model (\textbf{b}).
  }
 \label{err1}
\end{figure}

We now perform the proposed algorithm based on two different initialization of the low-fidelity model. The first analysis is based on the assumption that we already have some random dataset of solutions of a FOM. Then the initial low-fidelity model is a ROM which is constructed from this random sketch model.  In this example, we initialize the random sketch with 2 linearly independent snapshots from the training set. Since this is a low-rank linear problem, we chose to select one sample point per iteration. With a target tolerance set to $\epsilon= 10^{-6}$, 6 iteration cycles are required to achieve the desired accuracy, as shown in figure \ref{err1a}. Also, we show the convergence plot for parameters belonging to the validation set, which decays smoothly until the target accuracy is achieved. This implies, that the quality of ROM constructed with the proposed iterative multi-fidelity approach represents well the large-scale PDE system for any parameter belonging to the parametric space $\mathcal{D}$.

\begin{figure}[hbt!]
    \begin{subfigure}{0.4\textwidth}
    \includegraphics[trim=6cm 11cm 6cm 11cm, clip=true,width=\textwidth]{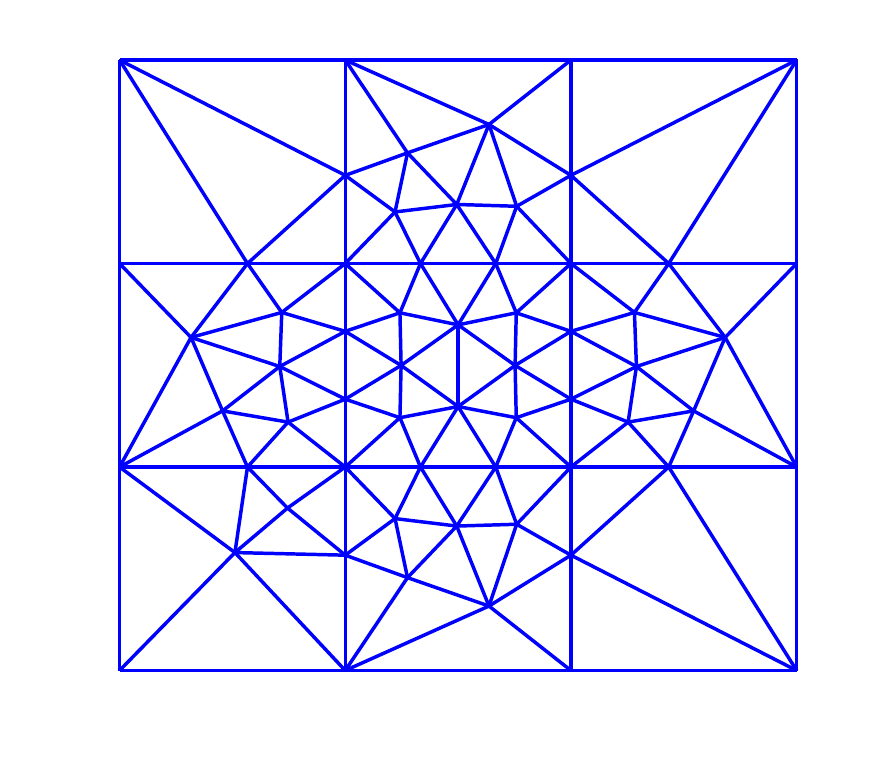}
       \subcaption{Coarse grid ($\#$ Nodes =62 )}
 \label{condmesha}

  \end{subfigure}
 \begin{subfigure}{0.4\textwidth}
   \includegraphics[trim=6cm 11cm 6cm 11cm, clip=true,width=\textwidth]{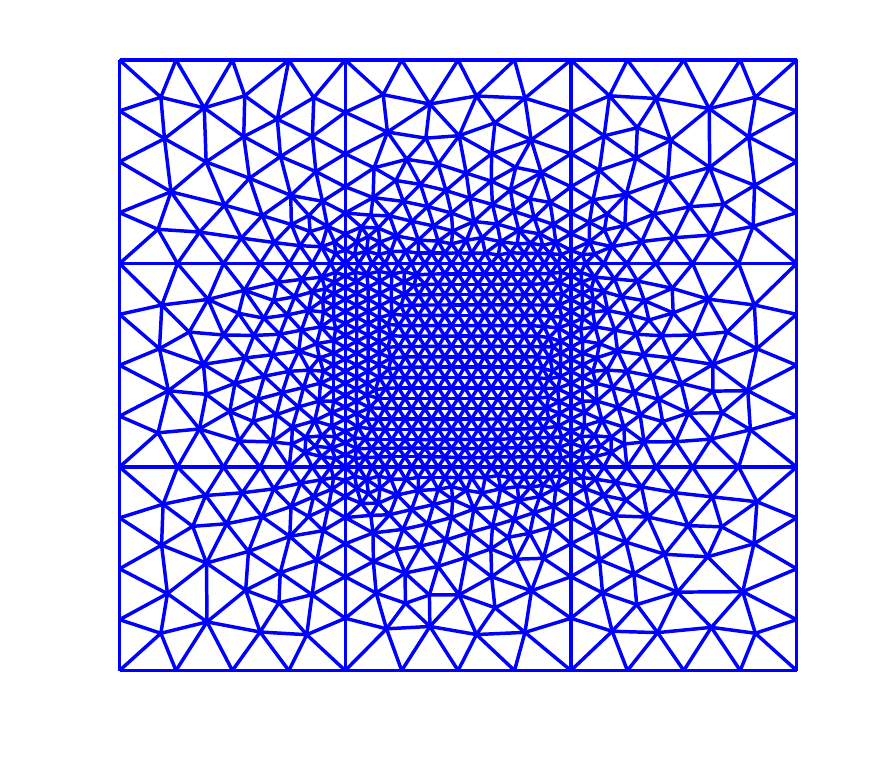}
  
        \subcaption{Fine grid ($\#$ Nodes =895)}
 \label{condmeshb}

  \end{subfigure}
      \caption{ (\textbf{a}) a Coarse discretized mesh for initial low-fidelity model, and (\textbf{b})  a fine discretized mesh for high-fidelity model approximations. } 
\end{figure}

 The second analysis is for the cases when no database of solutions is available a priori, we construct the initial low-fidelity model from a coarse sketch model as shown in figure \ref{condmesha}. Figure \ref{condmeshb} represents a fine grid model used for generating the high-fidelity solution. We observe that the same number of iterations are required as in the previous case to achieve accuracy of $\mathcal{O}(10^{-6})$ for parameters belonging to both the training and validation set, shown in figure \ref{err1b}. This specific example has only a parametric dimension of 2 and also this problem is not mesh-dependent, so both the sketch models consisted of sampling the same number of points. However, we will observe that the selection of points is not consistent in a complex problem as the advection-diffusion problem in 9 dimensions discussed in section \ref{numericaltest2}. Also, more or fewer sampling points may be required depending on the initialization of the low-fidelity model to retain the same target accuracy.
 
\begin{figure}[hbt!]
    \begin{subfigure}{0.4\textwidth}
     \includegraphics[trim=6cm 11cm 6cm 11cm, clip=true,width=\textwidth]{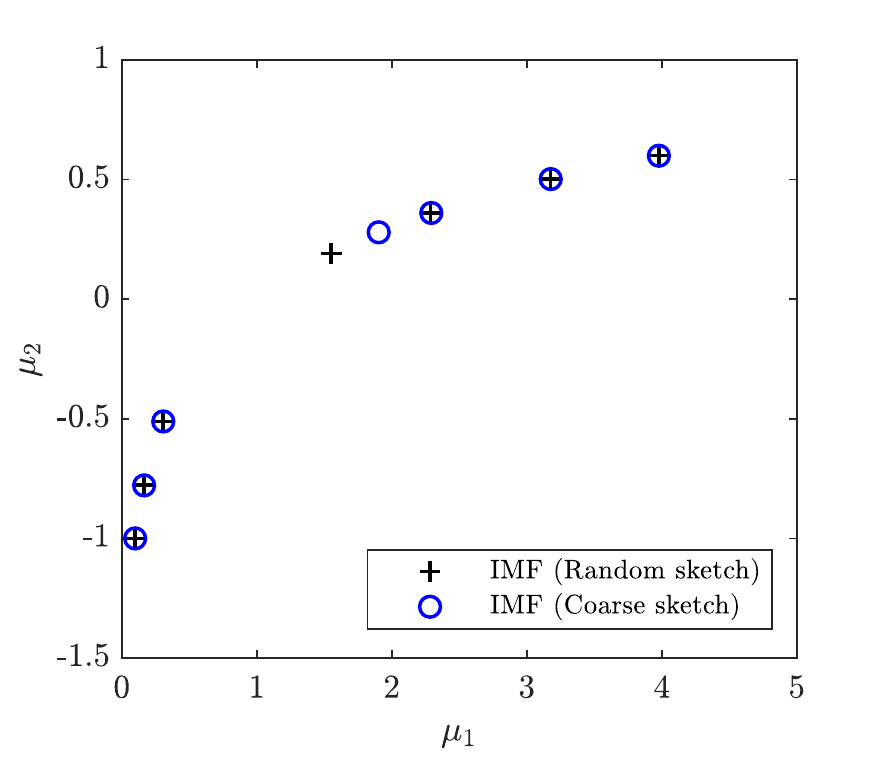}
     \subcaption{}
 \label{parametricptsa}
  \end{subfigure}
 \begin{subfigure}{0.4\textwidth}
   \includegraphics[trim=6cm 11cm 6cm 11cm, clip=true,width=\textwidth]{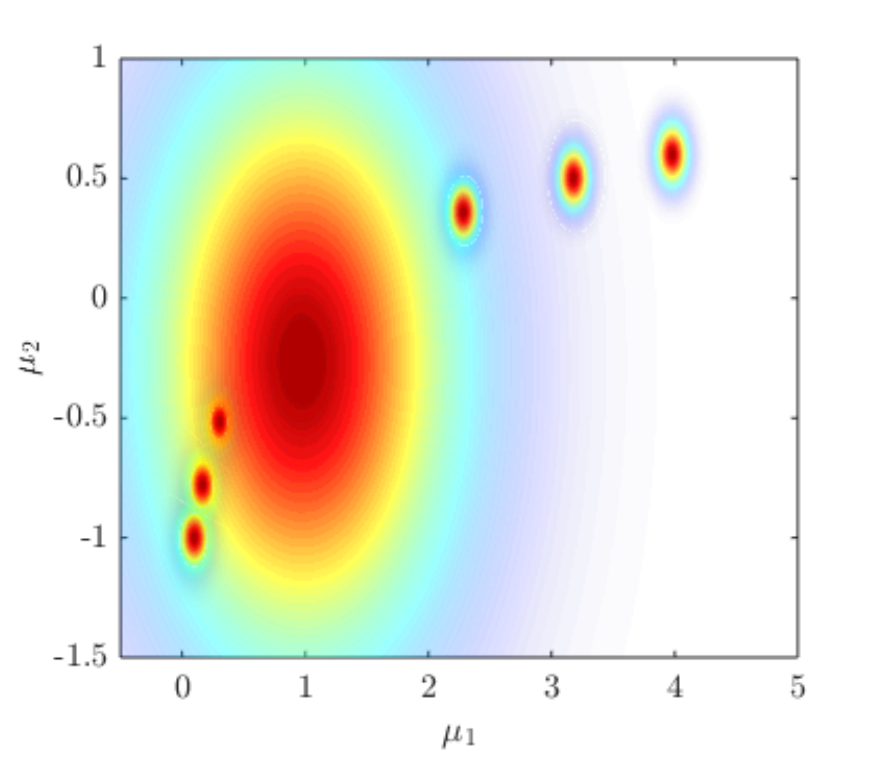}
   \subcaption{}
 \label{parametricptsb}
    
  \end{subfigure}
        \caption{ Sampled parametric points using two sketch models (\textbf{a}), Gaussian distribution of sampling points over 10 trials (\textbf{b}). }
  \end{figure}

   The parametric points sampled using both the sketch models are represented in figure \ref{parametricptsa}. We can notice that irrespective of the sketch model chosen to construct the low-fidelity model, the proposed method in this example extracted almost exact points in both cases. We also plot a density map in figure \ref{parametricptsb} showing the probability of a point to be picked at a certain location, which is obtained here by finding a Gaussian distribution over an ensemble of numerical experiments.  
   In other words, each experiment is initialized by different random snapshots without any repetition. The figure reflects a very interesting behavior, showing that six of the seven sampled locations are fairly the same in each trial except for one point that has more variance than the other six, which was noticed when the third point was chosen. This is attributed to the selection mechanism, since the third point is picked randomly after the first two points are drawn from the training set during the initialization of the random sketch model. Then accordingly the algorithm optimizes the location of the third point thus yielding consistent results. Since this approach is based on heuristics, obtaining exact sampling points in different cases is not of the utmost priority, as it is more on achieving unique points that can statistically well represent the entire parametric subspace.

  \begin{figure}[hbt!]
    \begin{subfigure}{0.4\textwidth}
   \includegraphics[trim=6cm 11cm 6cm 11cm, clip=true,width=\textwidth]{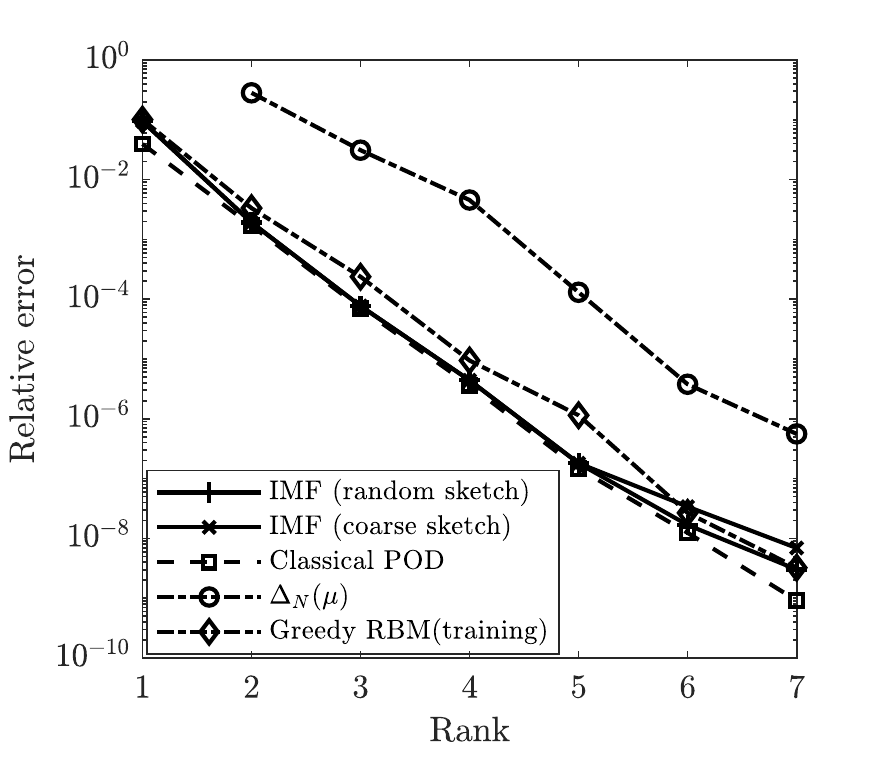}
        \subcaption{}
        \label{err2a}

 \end{subfigure}
        \begin{subfigure}{0.4\textwidth}
     \includegraphics[trim=6cm 11cm 6cm 11cm, clip=true,width=\textwidth]{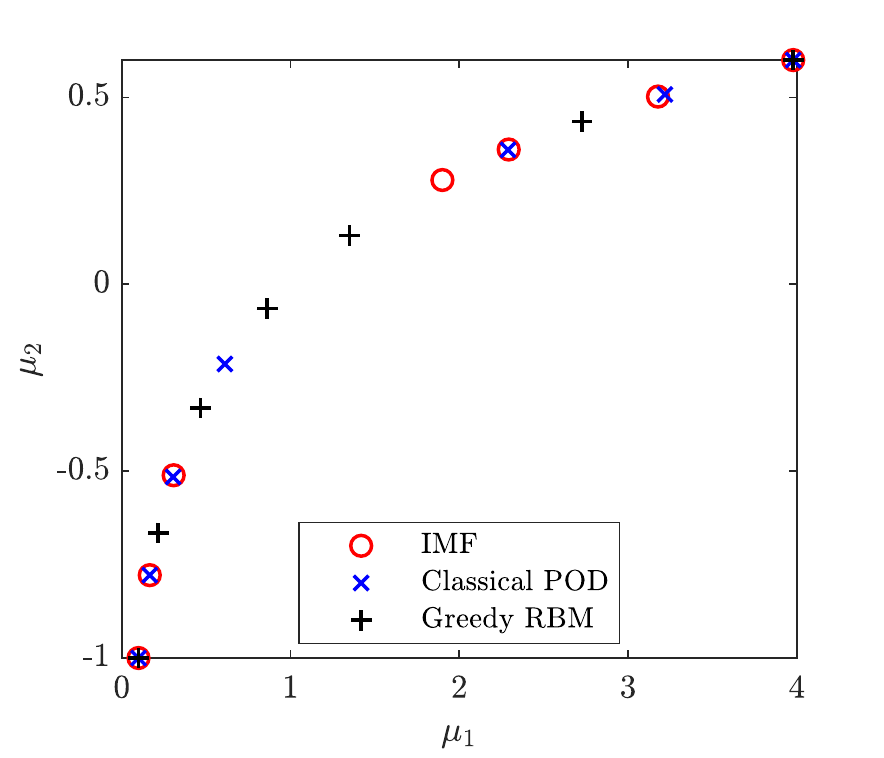}
        \subcaption{}
        \label{err2b}

  \end{subfigure}
  
\caption{    Error $\epsilon_{\text{\tiny POD}}$ of the POD projection, posteriori error bound by greedy RBM and error $\epsilon_{\text{\tiny ROM}}$ of the solution of the reduced order model obtained by the proposed method (\textbf{a}), and sampled parametric points using the proposed method, classical POD and greedy RBM (\textbf{b}).}
\end{figure}

We also plot the $\ell_2$ norm POD projection error ($\epsilon_{\text{\tiny POD}}$) and ROM error ($\epsilon_{\text{\tiny ROM}}$) against the rank of the low-fidelity model using both the sketch models in figure \ref{err2a}, as per the definitions in section \ref{errordef}. We can observe both the POD and ROM error  decay exponentially with an increase in the rank of the system. Also, the POD error curve is observed to be lower than the ROM error as expected, which is represented by the decay of singular values.  
Now we make a comparison by  solving the same problem using the greedy RBM algorithm, the lower bound is computed using the multi-min-theta approach explained in detail in \cite{Hesthaven2015}. The convergence plot for the relative $\ell_2$-norm ROM error  ($\epsilon_{\text{\tiny ROM}}$) for all the parameters in the training set and the error bound is shown in figure \ref{err2a}. It is observed that the ROM error has a smooth exponential decay and is lower than the max norm error as expected. 

For the target accuracy set to $\mathcal{O}(10^{-6})$, the rank of the system achieved is 7, implying that 7 parametric points were sampled to recover the basis functions similar to the findings of our proposed method. Thus, it is evident that the qualitative performance of the proposed method is comparable with the greedy RBM. We also show an illustration of sampled points obtained by greedy RBM  in figure \ref{err2b}. We notice that the parametric points picked in the proposed method are not entirely comparable to greedy RBM, however, the sampled points in both cases follow a logarithmic trend which may be related to the way we discretized our training set, where we chose to discretize  $\bm{\mu}_{[2]}$ using a log space and a uniform discretization for $\bm{\mu}_{[1]}$. If we apply DEIM on the parametric functions obtained by POD on high-fidelity snapshots for all $\bm{\mu} \in \Xi_{train}$ and plot the sampled points in the same figure \ref{err2b}, we observe that the selection of points is comparable with the proposed method and lie on the same logarithmic trend.

 \section{9D advection-diffusion problem }
 \label{numericaltest2}
\begin{figure}[hbt!]
      \includegraphics[trim=0cm 0cm 0cm 0cm, clip=true,width=0.8\textwidth]{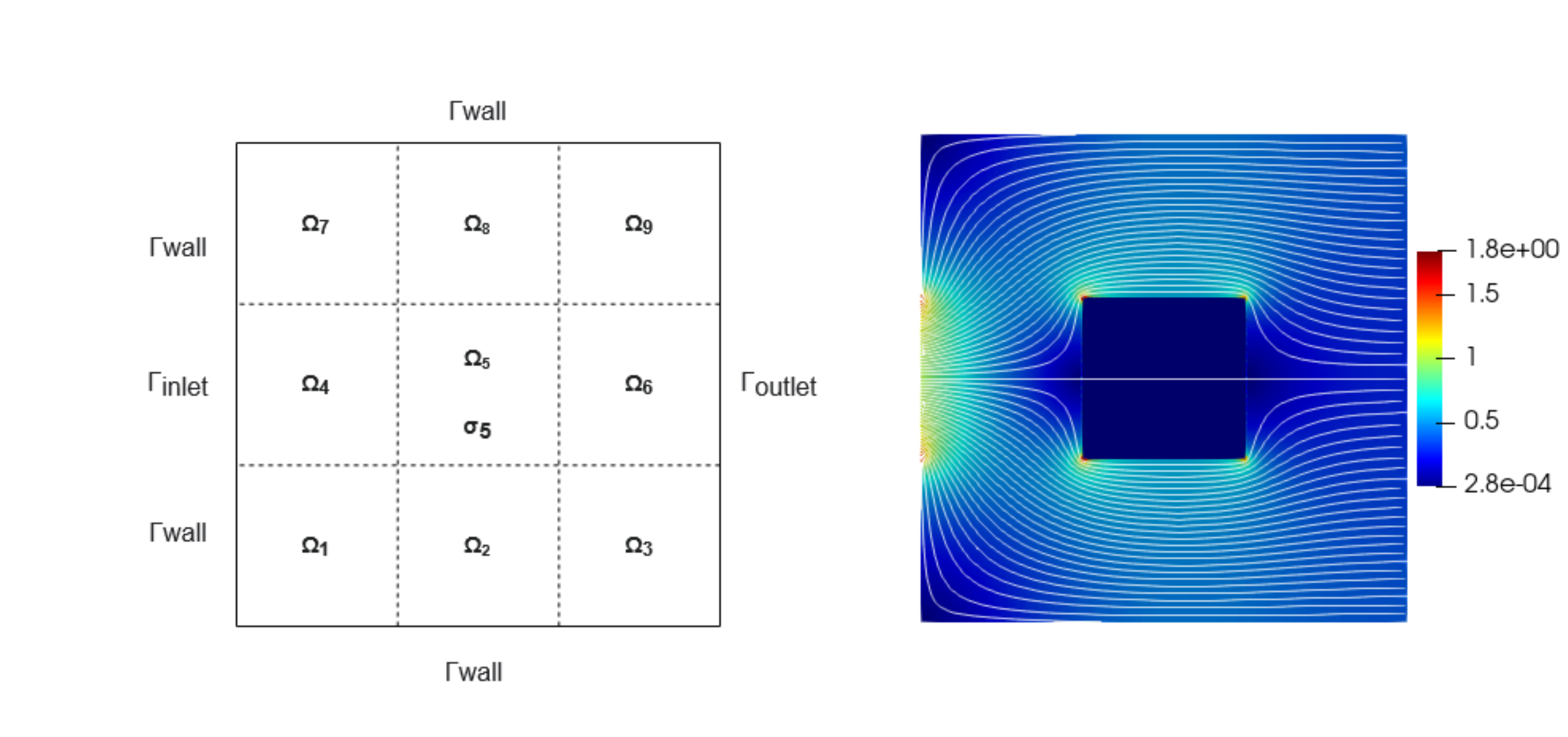}
       \caption{ (left) Geometrical set up of advection-diffusion problem in a 9 block system and (right) advective flow field }
\label{CVR}
\end{figure}

In this section, we study an advection-diffusion problem with a source term for a 9 parametric dimension. As it can be seen in figure \ref{CVR}, the domain is divided into 9 subdomains where each region has a different diffusivity coefficient which serves as the input parameter.

\subsection{Problem setting}
The strong form of the parameterized advection-diffusion reaction equation is governed by the elliptic PDE. For some parameter value, $\bm{\mu} \in \mathcal{D}$ find $u(\bm{\mu})$ such that:

\begin{equation}   
\begin{aligned}
\bm{b}\cdot \nabla u(\bm{\mu})-\mathbb{K}(\mu_{i}) \Delta u(\bm{\mu})+\sigma & =0 \qquad  \forall\Omega_{i} \quad i=1, \ldots, 9 \\
u(\bm{\mu})     & = 0  \qquad \text{ at } \Gamma_{inlet}\\
\nabla u(\bm{\mu}) \cdot \hat{n} & = 0 \qquad \text{ at } \Gamma_{walls}\\
\sigma & = 1 \qquad \text{ at } \Omega_5\\
\label{eq:cdreq}   
\end{aligned}
\end{equation}
        where, $u(\bm{\mu})$ is the unknown field variable  for a spatial domain $\Omega=[0,1]\times [0,1]$. The diffusion coefficient is given by
    $\mathbb{K}(\bm{\mu})=\bm{\mu}$  for the input  parameter $\bm{\mu} \in \mathcal{D} = [0.01,10]^9$ and $\sigma$ is the constant reaction term.
    $\bm{b}$ is a given-advected flow field which is  obtained by solving a potential flow problem in the same bounded domain such that $\nabla. \bm{ b}=0$. We consider a low permeability zone in $\Omega_5$, and relatively higher permeability over the rest of the domains, and with a constant velocity at the inlet, we obtain the streamlines as shown in figure \ref{CVR}. 
    
The weak parameterized formulation using SUPG then reads as: for some parameter $\bm{\mu} \in \mathcal{D}$, find $u(\bm{x};\bm{\mu}) \in \mathcal{V}^\mathcal{N}$ where $\mathcal{V}^\mathcal{N}=\left\{v \in (H_0^1(\Omega))^2|v|_{\Gamma_{\text {inlet}}}=0\right\}$,

\begin{equation}
 \langle v+\frac{\beta h}{2} \nabla v, \bm{b} \cdot \nabla u\rangle+\langle\nabla v, \mathbb{K}(\bm{\mu}) \nabla u\rangle -\langle \frac{\beta h}{2} \nabla v, \mathbb{K}(\bm{\mu}) \Delta u\rangle 
 +\langle  v+\frac{\beta h}{2} \nabla v, \sigma \rangle =0 \qquad   \forall v \in \mathcal{V}^\mathcal{N}
  \label{eq:cdrweak}
\end{equation}
 The parameter $\beta$ is a dimensionless constant that depends on the so-called Peclet number given by $Pe=\frac{||\bm{b}||h}{2\mathbb{K}(\bm{\mu})}$, where $h$ is a typical element size in the direction of the velocity and $||\bm{b}||$ is the characteristic flow velocity. Note that the third term in the equation \ref{eq:cdrweak} is zero in the case of linear elements used for domain discretization \cite{Hulsen2009}.

The reduced weak form of eq. \ref{eq:cdrweak} is obtained by projection onto the properly selected low-dimensional subspace spanned by a reduced basis function $\{\bm{\phi}^k\}_{k=1}^{r}$ such that,  

\begin{equation}
 \langle \phi^k+\frac{\beta h}{2} \nabla \phi^k, \bm{b} \cdot \nabla u_{\text{\tiny POD}} \rangle+\langle\nabla \phi^k, \mathbb{K}(\bm{\mu}) \nabla u_{\text{\tiny POD}}\rangle 
 +\langle  \phi^k+\frac{\beta h}{2} \nabla \phi^k, \sigma \rangle =0 \qquad    \forall k=1,\ldots,r
  \label{eq:cdrROM}
\end{equation}

where, $$u_{\text{\tiny POD}}(\bm{x};\bm{\mu})=\sum_{j=1}^r \phi^{j}(\bm{x})z_{j}(\bm{\mu}) $$

$ \bm{z}=\{z_1, \ldots, z_r\}^{\text{\tiny T}} $  represents the  coefficients of the POD expansion.

Note that only the diffusion term in the  equation \ref{eq:cdrROM} is affine with respect to the input parameter $\mathbb{K}(\bm{\mu})$  and can be efficiently reduced during the offline stage,
 
 \begin{equation}
 \begin{gathered}
 \langle \nabla \phi^k,\mathbb{K}(\bm{\mu}) \nabla u_{\text{\tiny POD}}\rangle = \mathbb{K}(\bm{\mu})\sum_{j=1}^{r}  \langle \nabla \phi^k, \nabla \phi^j\rangle z_j \, ,
  \label{eq:exred}
\end{gathered}
\end{equation}
 
  in which the  $(r \times r)$ operator $\langle \nabla \phi^k,\mathbb{K}(\bm{\mu}) \nabla u_{\text{\tiny POD}}\rangle$ can be computed once and for all in the offline stage. During the online stage, if new parameter $\mathbb{K}(\bm{\mu})$ is prescribed, the evaluation of the diffusion operator of eq. \ref{eq:cdrROM} can be done in reduced complexity (i.e. it does not depend on the original dimension $\mathcal{N}$) since it only requires $\mathcal{O}(r \times r)$ operations. This step is crucial for retaining the computational efficiency of the ROM, however, the same idea cannot be straightforwardly applied to the convective and source term of the equation \ref{eq:cdrROM} 
  as the projection operator dependency on the input parameter i.e. diffusion coefficient is non-affine. However, the non-affineness is not addressed in this study; instead, the application of the suggested approach to a high-dimensional PDE system is the main focus.  

\subsection{Results and discussion}
\begin{figure}[hbt!]
   \begin{subfigure}{0.3\textwidth}
     \includegraphics[trim=12cm 7cm 12cm 5cm, clip=true,width=\textwidth]{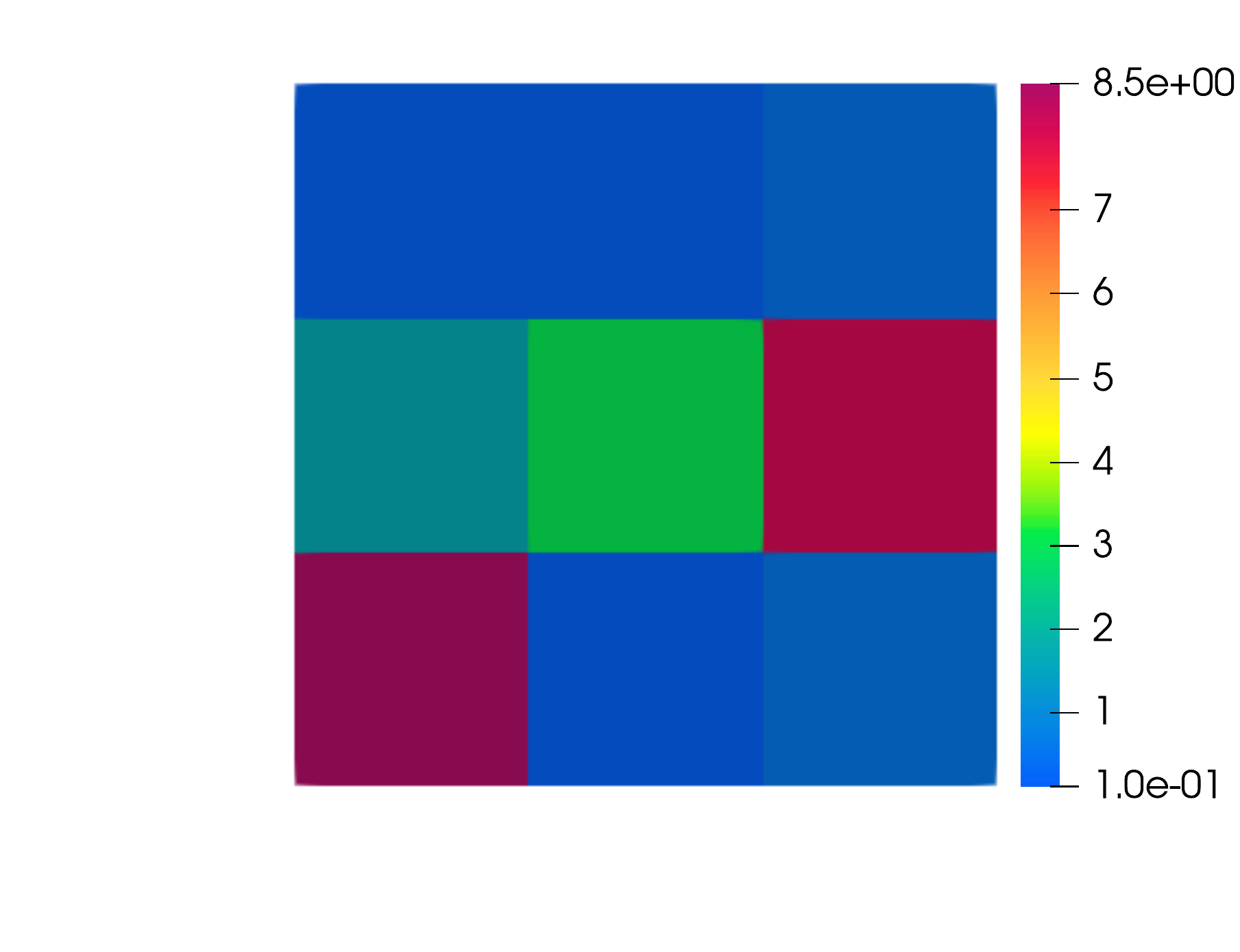}

         \end{subfigure}
    \begin{subfigure}{0.3\textwidth}
          \includegraphics[trim=12cm 7cm 12cm 5cm, clip=true,width=\textwidth]{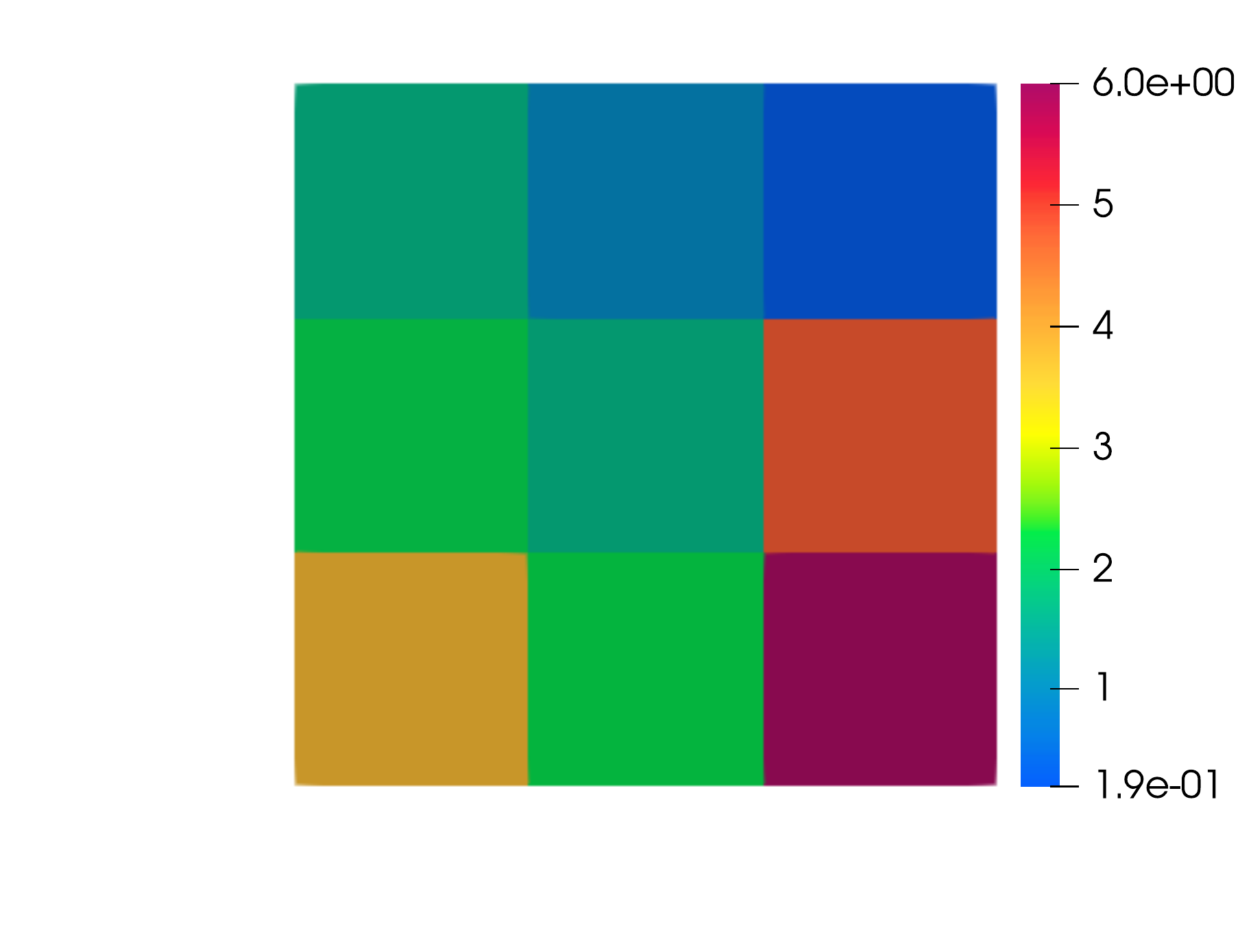}
    \end{subfigure}
      \begin{subfigure}{0.3\textwidth}
          \includegraphics[trim=12cm 7cm 12cm 5cm, clip=true,width=\textwidth]{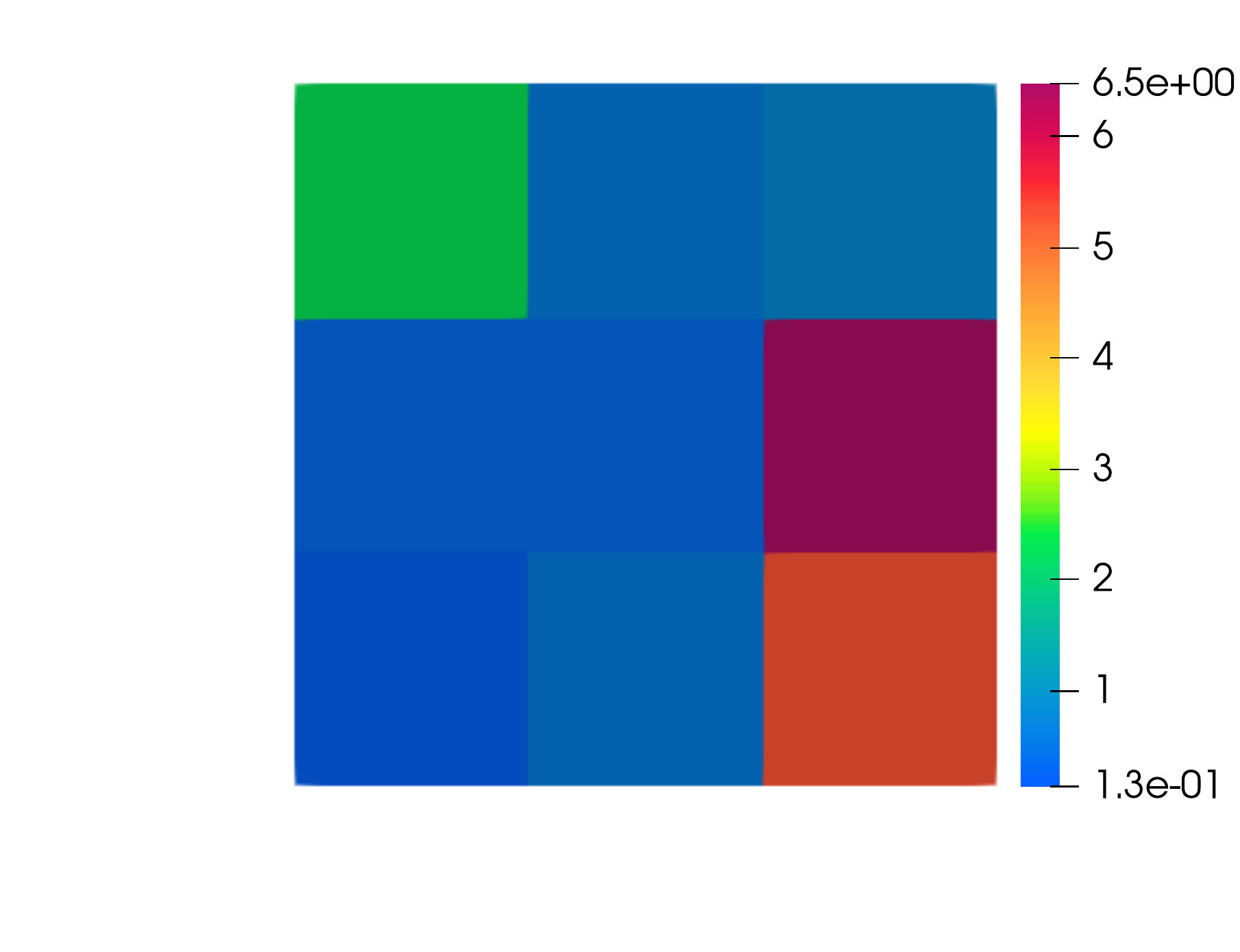}
    \end{subfigure}
    
   \begin{subfigure}{0.3\textwidth}
     \includegraphics[trim=12cm 7cm 12cm 4cm, clip=true,width=\textwidth]{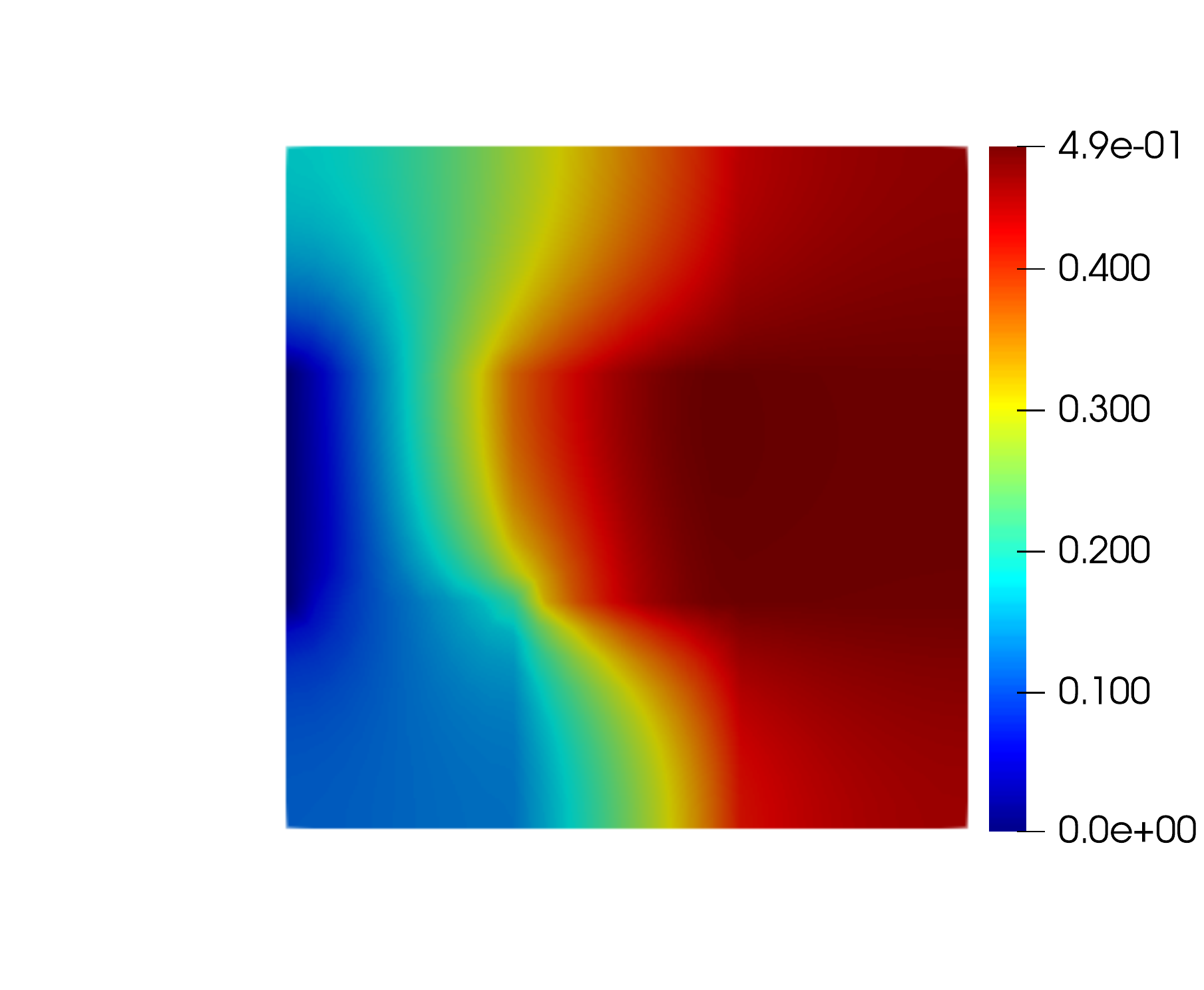}
         \end{subfigure}
    \begin{subfigure}{0.3\textwidth}
          \includegraphics[trim=12cm 7cm 12cm 4cm, clip=true,width=\textwidth]{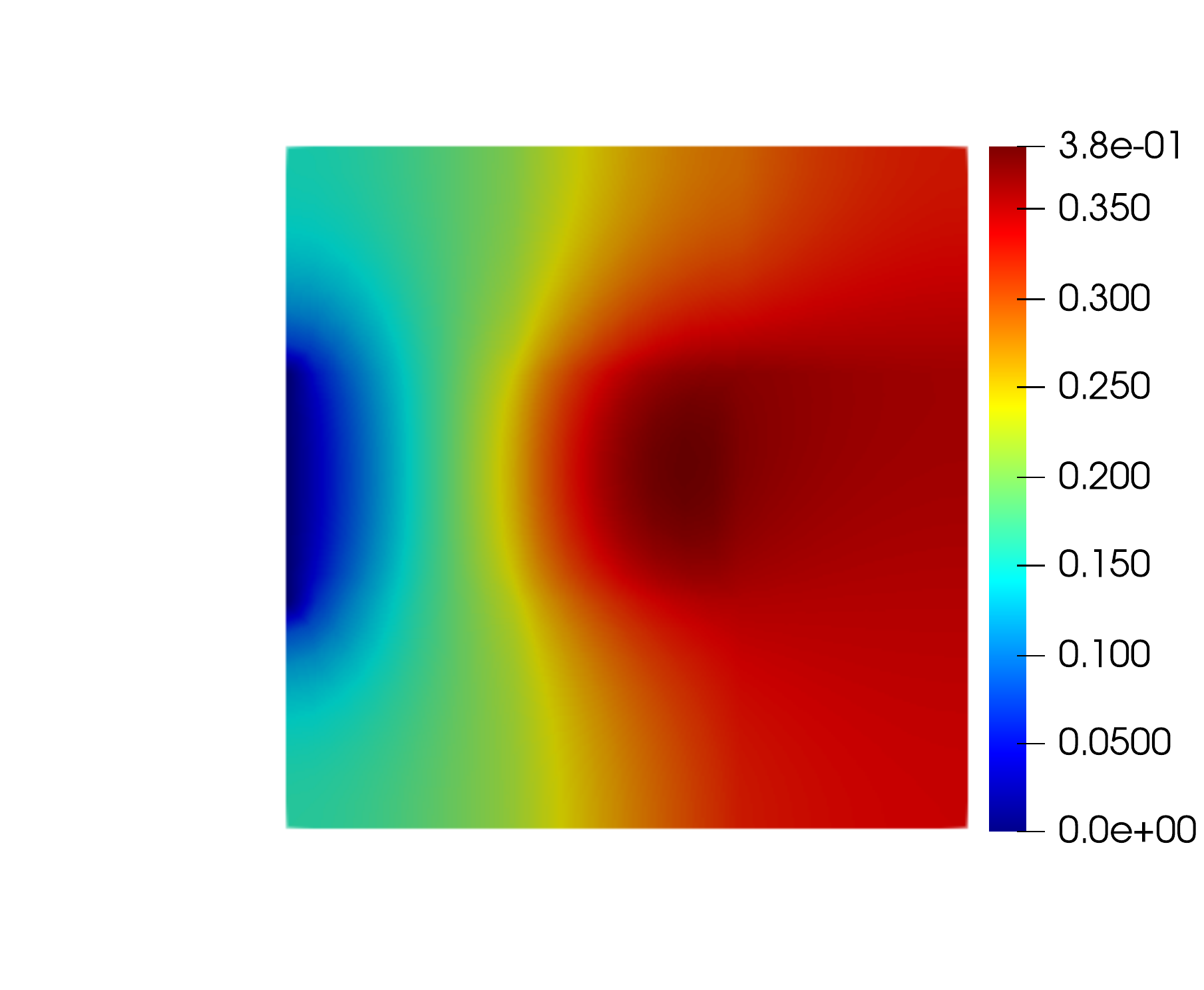}
    \end{subfigure}
      \begin{subfigure}{0.3\textwidth}
          \includegraphics[trim=12cm 7cm 12cm 4cm, clip=true,width=\textwidth]{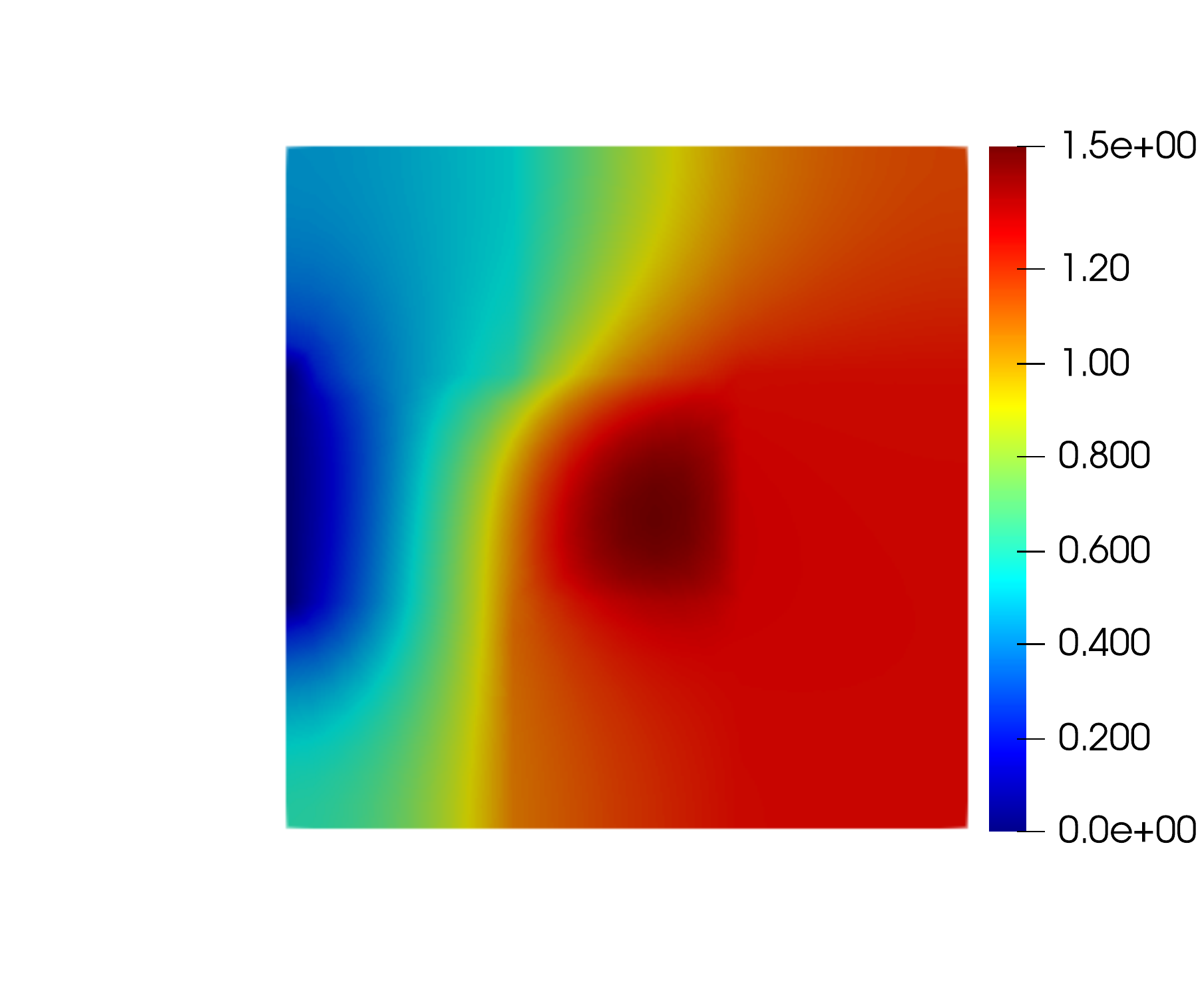}
    \end{subfigure}

    \caption{Three different representative solutions for the advection-diffusion  problem shown (below) by varying nine different combinations of diffusion coefficient $\mathbb{K}(\mu)$ in all the three domains (top).} 
  \label{transport1}
\end{figure}

The input parameter $\bm{\mu}$ is  discretized using the LHS technique with 2500 sample points, from which the training set $\Xi_{train} \subset \mathcal{D}$ consist of 2000 points and the remaining 500 samples are used for the validation set $\Xi_{val} \subset \mathcal{D}$ to certify the quality of reduced basis approximation. Figure \ref{transport1} represents the solution field with different combinations of diffusion coefficients for each of the 9 blocks.

\begin{figure}[h!]
\centering
  \begin{subfigure}{0.4\textwidth}
    \includegraphics[trim=6cm 11cm 6cm 10cm, clip=true,width=\textwidth]{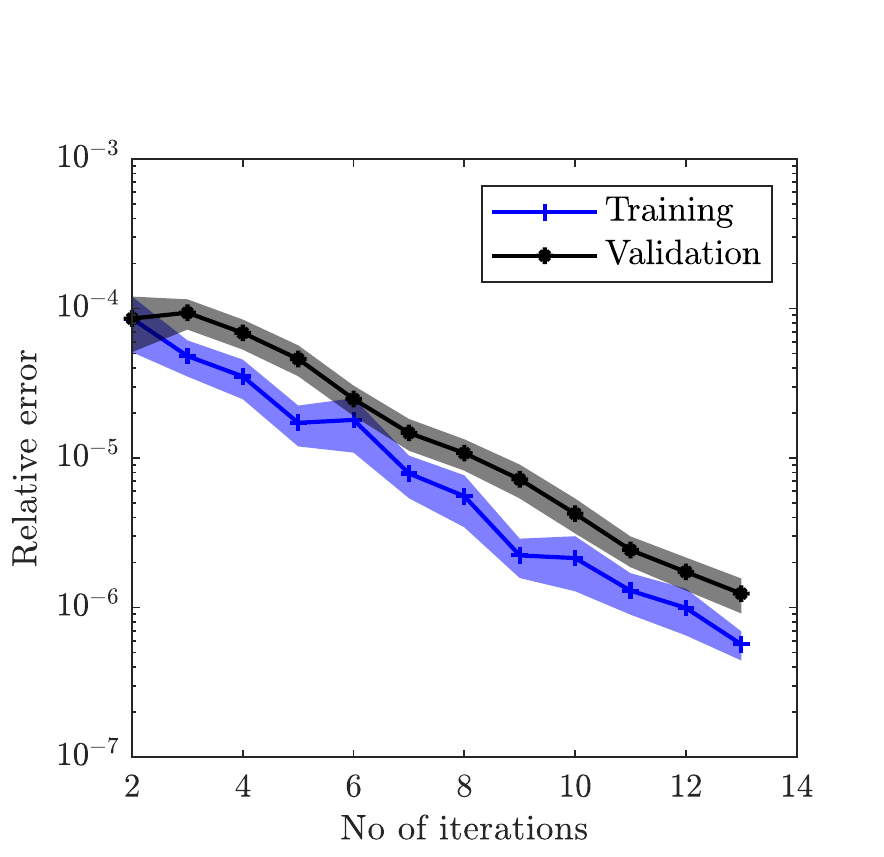}
    \subcaption{\tiny initial rank 100, incremented by 10 points every iteration }
     \label{CDR1a}
  \end{subfigure}
    \begin{subfigure}{0.4\textwidth}
         \includegraphics[trim=6cm 11cm 6cm 10cm, clip=true,width=\textwidth]{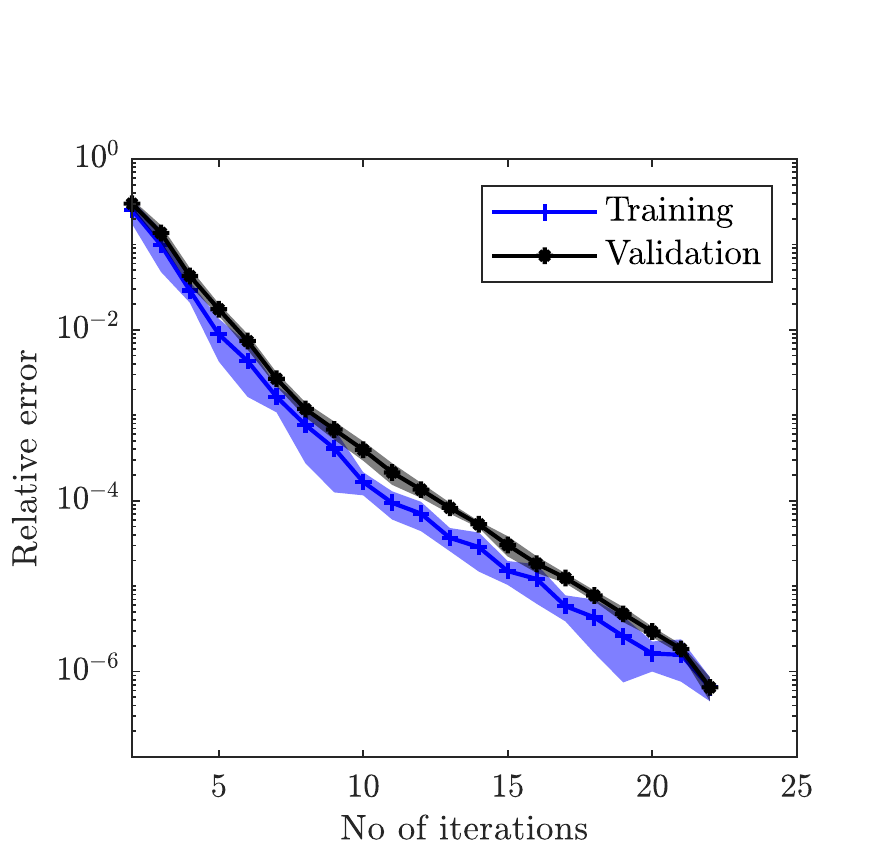}
      \subcaption{\tiny initial rank 10, incremented by 10 points every iteration}
       \label{CDR1b}
  \end{subfigure}
  \hfill
     \begin{subfigure}[b]{0.4\textwidth}
\centering
 \includegraphics[trim=6cm 11cm 6cm 10cm, clip=true,width=\textwidth]{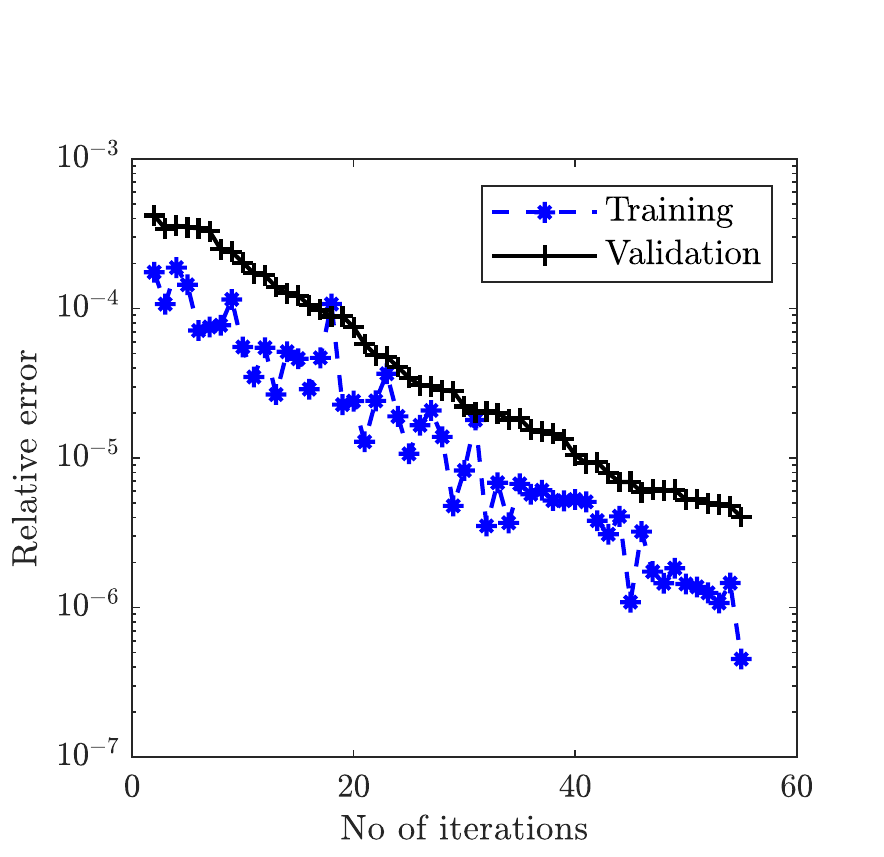}
       \subcaption{\tiny initial rank 100, incremented by 2 points every iteration}
 \label{CDR1c}

\end{subfigure}

     \caption{(\textbf{a,b}) Mean of error $\epsilon_{train}$ and $\epsilon_{val}$ between FOM and ROM solution obtained using different random sketch models over 10 trials and its std. deviation. (\textbf{c}) Error  $\epsilon_{train}$ and $\epsilon_{val}$ between FOM and ROM solution for one trial. }  
 \label{CDR1}
  \end{figure}

 Similar to the previous numerical example, we perform the algorithm based on two different initialization of the low-fidelity model. First, the discussion is presented for the low-fidelity model approximated using a random sketch model. Three studies are conducted to evaluate the computational performance: the first two examined the impact of the random sketch on the qualitative performance of the suggested method, while the third examined the quantitative performance of the algorithm by adjusting the greedy parameter for sampling. This algorithm is conducted over 10 trials for the first two cases, the convergence plot is shown in figure \ref{CDR1}. In the first case study, the initial rank of the random sketch for constructing a ROM is chosen to be 100, and during every iteration cycle, 10 points are added incrementally to update the low-fidelity model.  With a target tolerance set to $\epsilon= 10^{-6}$, we can observe that 13 iteration cycles are required to achieve the desired accuracy in each of the 10 trials, refer to figure \ref{CDR1a}. A total count of $100+(10*12) = 220$  parametric points are sampled out of 2000 points from the training set. It is to be noted that the rank of the low-fidelity model is also enhanced by 10 which implies that all the points sampled are unique which is expected and as a consequence, the recovered basis functions by construction are linearly independent. In the second case study, the random sketch is initialized with 10 linearly independent snapshots from the training set instead of 100 and in each iteration, 10 points are sequentially added to recover the basis functions. It is observed from figure \ref{CDR1b}, that in 22 iterations the target accuracy is achieved with a total sampling of  220 points out of 2000 points from the training parametric set, similar to the first case. Thus, it is evident that irrespective of the size of the random sketch chosen for the initial construction of the low-fidelity model, the algorithm performed well in both scenarios and the final enrichment of the low-fidelity model converged towards the FOM model accurately within the prescribed tolerance. The computational time required to achieve the target accuracy  is of the same order  $\mathcal{O}(10^3)$ in seconds in both cases, while having comparable computational performance.
 
In the third case, the study is conducted by sampling 2 parametric points per iteration instead of 10 points.  Figure \ref{CDR1c} reflects 54 iteration cycles required to achieve the same target accuracy, with a total sampling of  (100+54*2=208) points, unlike in the previous two cases where 220 points were selected from the parametric space. The CPU time required in this case is $\mathcal{O}(10^4)$ seconds, which is one order higher than the previous cases. This implies adding a few points per iteration can minimize the risk of sampling excess points while maintaining the same order of accuracy, but at the cost of higher CPU time. Due to the discrete nature of error evaluation, the relative training error is observed to be noisy, but with the validation error plot, we can see a smooth decay of the curve as the error is evaluated over the entire validation set rather than at select discrete points.

\begin{figure}[hbt!] 
  \begin{subfigure}{0.4\textwidth}
     \includegraphics[trim=6cm 10cm 6cm 10cm, clip=true,width=\textwidth]{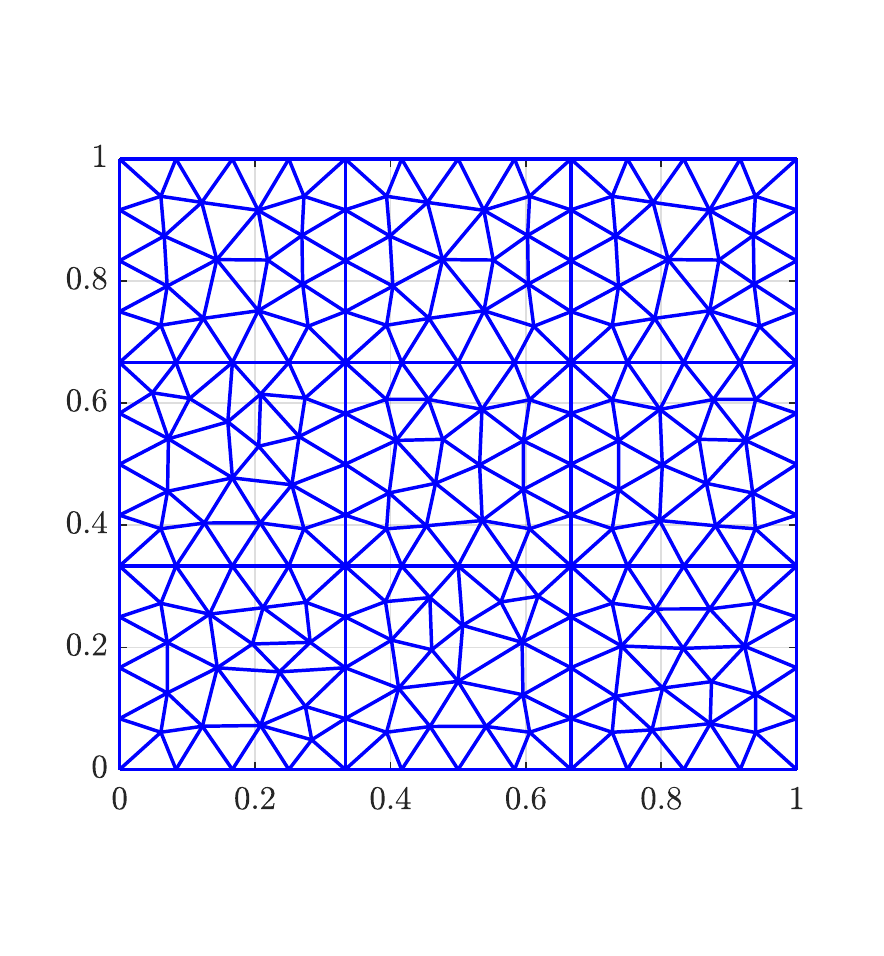}
     \subcaption{Coarse grid ($\#$ Nodes=575)}
     \label{HTmesha}
    \end{subfigure}
  \begin{subfigure}{0.4\textwidth}
 \includegraphics[trim=6cm 10cm 6cm 10cm, clip=true,width=\textwidth]{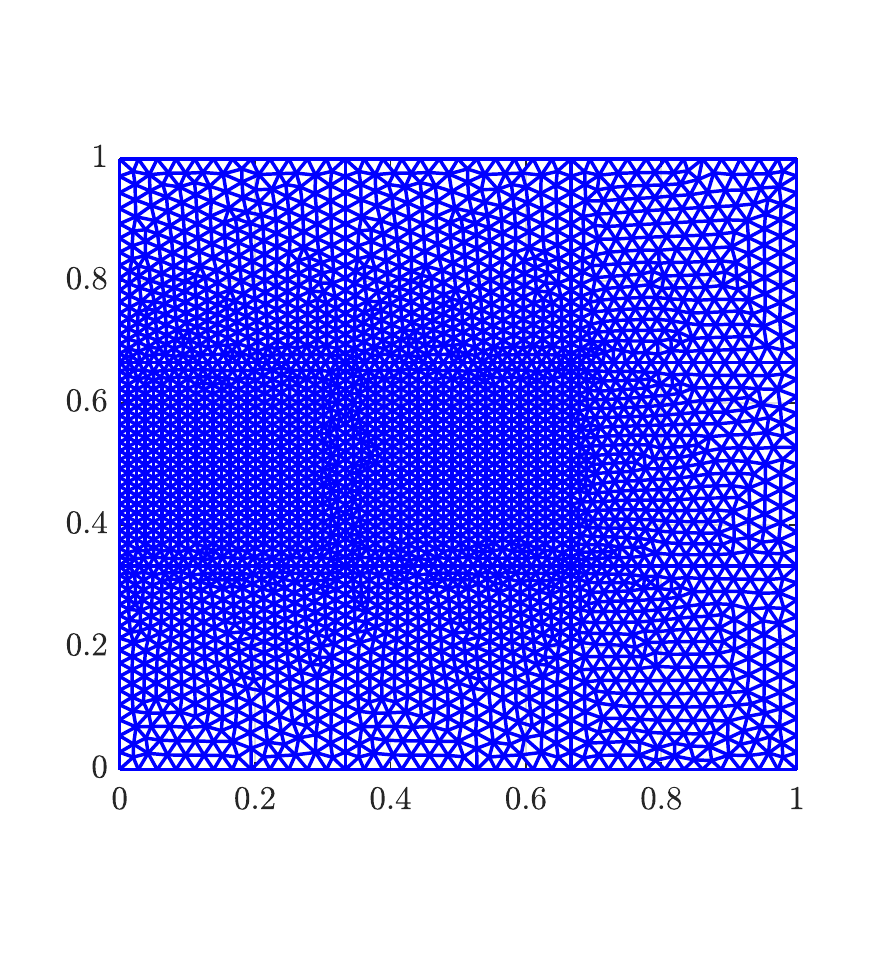}
 \subcaption{Fine grid ($\#$ Nodes=3492)}
 \label{HTmeshb}
  \end{subfigure}
  \caption{(\textbf{a}) Coarse discretized mesh for initial low-fidelity model, and (\textbf{b}) fine discretized mesh for high-fidelity model  approximations.}
\end{figure}

\begin{figure}[hbt!]
  \begin{subfigure}{0.41\textwidth}
     \includegraphics[trim=6cm 10cm 6cm 10cm, clip=true,width=\textwidth]{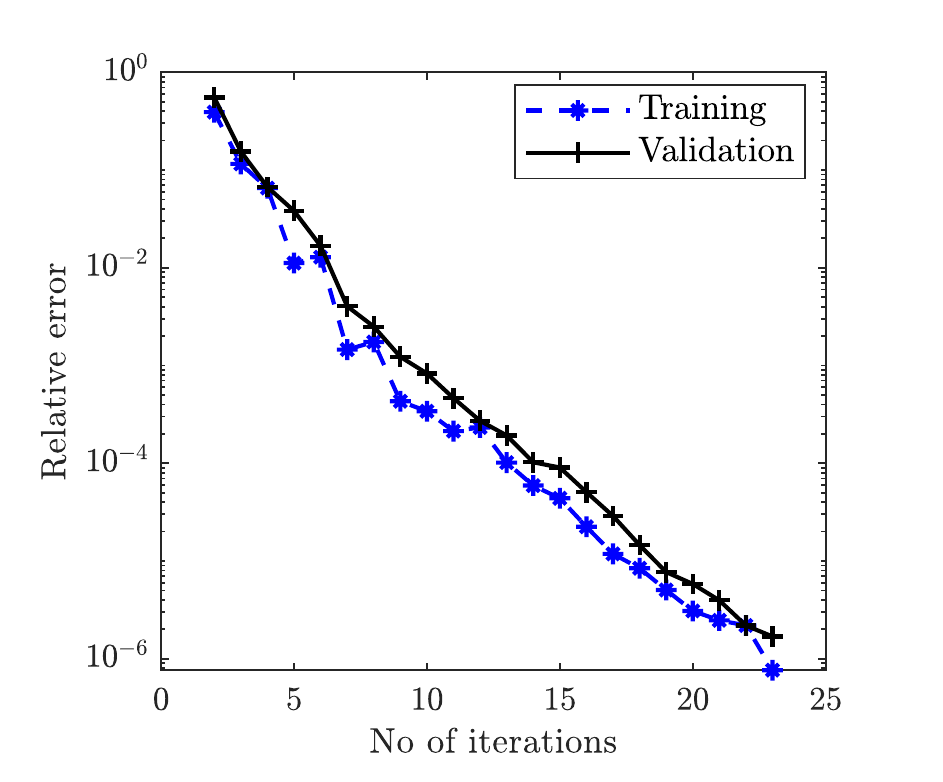}
     \subcaption{}
      \label{CSMHTa}

  \end{subfigure}
      \begin{subfigure}{0.39\textwidth}
     \includegraphics[trim=6cm 11cm 6cm 11cm, clip=true,width=\textwidth]{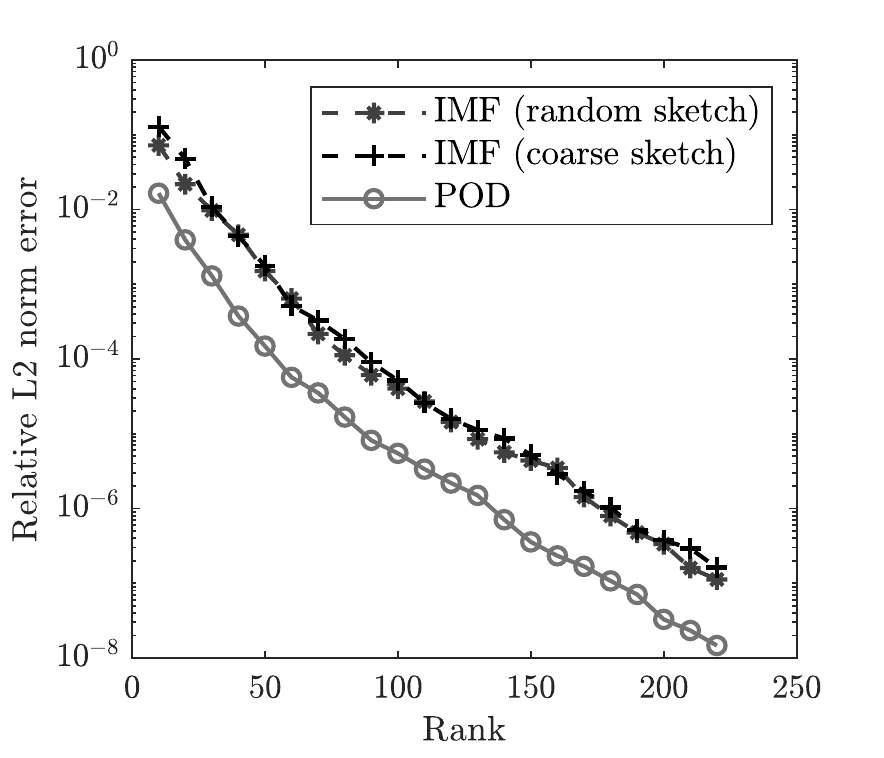}
     \subcaption{}
      \label{CSMHTb}

  \end{subfigure}
  \caption{ Error $\epsilon_{train}$ and $\epsilon_{val}$ between FOM and ROM solution using coarse sketch model (\textbf{a}).
     Error $\epsilon_{\text{\tiny POD}}$ of the POD projection, and error $\epsilon_{\text{\tiny ROM}}$ of the solution of the reduced order model obtained by the proposed method \textbf{ (b)}.}
  \end{figure}

      For the second analysis, the initial low-fidelity model is built using a very coarse sketch model, as shown in figure \ref{HTmesha}. For recovering the high-fidelity solution, a fine discretized model is used in figure \ref{HTmeshb}. The points are added sequentially by incrementing with 10 every iteration. The target accuracy is achieved in 23 iterations (total count of sampled points is 230) as shown in figure \ref{CSMHTa}, which is more by 10 points compared to the random sketch model. As is already discussed, such types of PDE problems face numerical stability issues in case of high Peclet number (advection-dominated cases) and can be resolved  by applying artificial diffusion in the upwind direction. It is also important to note that the amount of artificial diffusion added to the system depends on the mesh size. So coarser the mesh size, the more diffusion is required, which may affect the solution significantly. Hence, to accurately capture the physical properties of the PDE system, it is very important to consider that the initial grid takes into account all the physical aspects of the problem such that there is no loss of information. This could be the plausible explanation for a higher number of points needed when a coarse sketch model is used as compared to the ROM for the initial construction of the low-fidelity model. Nevertheless, our main objective to guarantee convergence is achieved irrespective of the initial size of the coarse mesh used and can be confirmed further from the decay of the validation error curve seen in figure \ref{CSMHTb}. The CPU time taken is of the $\mathcal{O}(10^3)$ seconds, same as the first and second case. 
      
   Figure \ref{CSMHTb} shows the plot for $\ell_2$ norm POD projection error and ROM error for all the parameters belonging to the training set using both sketch models.   Both the POD and ROM errors have an exponential decay as the rank of the system increases, with the POD error serving as a lower bound to the ROM error. This proves the reliability of the proposed method on the quality of ROM constructed such that irrespective of the initial design  of the low-fidelity model, the ROM error displays similar decay properties as the POD error.

\section{Conclusions}

In this work, we presented the feasibility of a multi-fidelity approach in reduced basis approximation for a multi-dimensional parametric PDE system in an iterative procedure. The parametric sampling is adhoc and extracted from  low-fidelity model approximations based on heuristic assumptions. Such approximations of low-accuracy low-fidelity models over high-fidelity models enhance the computational performance in the offline stage significantly. Not to mention, this approach proved to be adequate when a posteriori error estimator is unavailable, which is an essential ingredient for efficient greedy sampling. The greedy selection is user dependent, i.e. one can select a single point or multiple points for each iteration. However, attention is to be given as excess parametric points might be sampled than required to achieve the same target accuracy but at a higher CPU cost. So a compromise can be made on the trade-off between computational efficiency and accuracy.

Overall, the construction of reduced basis subspace required less high-fidelity snapshot generation in the proposed method. This methodology is successfully demonstrated on a 2D steady-state advection-diffusion problem for 9 input parametric dimensions. A qualitative comparison is also presented for a simple steady-state heat conduction problem between the proposed method and greedy RBM, in which both cases had comparable computational performance. We also presented two different ways of initializing a low-fidelity model and irrespective of the initial quality of the low-fidelity model approximation, the method is shown to be reliable and stable by converging towards the FOM approximation within the prescribed tolerance. 

In our current work, we have only performed linear reduction, but we can further embed hyper reduction in our current framework "on the fly" i.e. during the multi-fidelity iterations for treating non-affine problems. This step can alleviate the cost of low-fidelity model approximations significantly. For some of the current research work one can refer to the articles \cite{alla2022adaptive,wen2022globally} on adaptive hyper reduction techniques which allows enrichment of the reduced integration domain during the online stage as the simulation progresses.

Additionally, during the multi-fidelity iterations, we may also evaluate the low-fidelity model solely on a portion of the randomly chosen parametric points as opposed to the complete training set. This procedure could drastically improve the computational performance of the methodology, specially for nonlinear PDEs.  However, there are certain implications to it such as there will be missing information in the parametric subspace and as suggested in the methodology, to sample points using the DEIM strategy, one needs to orthogonalize the current parametric functions in relation to the previous ones. As a possible solution to this problem, Gappy-POD may be used to reconstruct the missing data in the updated parametric functions and can be implemented into this technique for future research work.

\newpage

\bibliographystyle{unsrt}
\bibliography{ref}

\end{document}

%% file: sequence.tex
\coordinate (b) at (0,0);
\coordinate (c) at ($(b)+(-24,0)$);
\coordinate (d) at ($(b)+(20,-12)$);
\coordinate (e) at ($(b)+(20,-30)$);
\coordinate (f) at ($(b)+(-16,-40)$);
\coordinate (g) at ($(b)+(-16,-28)$);
\coordinate (h) at ($(b)+(-16,-12)$);
\coordinate (i) at ($(b)+(-30,-12)$);

\coordinate (o) at ($(b)+(0,-20)$);


\node[ellipse,draw,align=center, minimum width=1cm,minimum height=4cm,fill=darkgray!20](step1) at (c){ Training set,\\ Validation set, \\ Lo-Fi model  \\
Hi-Fi model};


\node[blueblock, align=center] (step2)at (b){
           \textbf{Step 1} \\ Solve the initial lo-fi model on training set
    };

      \node[rectangle,
    draw = black,
    dashed,
    text = black,
    fill = gray,
    minimum width = 70cm,
    minimum height = 50cm, fill opacity = 0.1] (rect1) at (o) {};

   \node[text width=4cm](phi) at (i) {STOP};

\node[blueblock](step3)at(d){\textbf{Step 2} \\ Estimate parametric modes};
    
\node[blueblock](step4)at(e){\textbf{Step 3} \\ Select new point(s) in the parametric space using DEIM on parametric modes};

\node[blueblock](step5)at(f){
\textbf{Step 4}\\
Solve hi-fi model on selected points and enrich the model's basis functions

};

\node[blueblock](step6)at(g){
\textbf{Step 5}\\ 
Solve a physics based reduced problem on the training set
}; 
\node[decision, minimum width=3cm, align=center,fill=red!10](step7) at (h){\textbf{Step 6}\\ 
Check training \\ and validation error \\ on select points};

\draw[myarrows]  (step1.east) -- (step2.west);
\draw[myarrows]  (step2.east) -|(step3.north);

\draw[myarrows] (step3.south) -- (step4.north);

\draw[myarrows]  ( step4.south) |-(step5.east);
\draw[myarrows]  (step5.north) --(step6.south);

\draw[myarrows]  (step6.north) --(step7.south);

\draw[myarrows]  (step7.east)  --node[xshift=-4cm,above=0.2cm]{ $< $ tol}node[xshift=-4cm,below=0.2cm,align=center]{ Update lo-fi \\ with ROM}(step3.west);

\draw [myarrows] (step7.west)--node[above=0.2cm,align=center]{ $>$ tol}(phi);